\documentclass[10pt]{article}

\usepackage[dvips]{epsfig} % this package is for manipulating postscript graphics
\usepackage{amssymb,amsmath,amsthm,mathrsfs} %this package contains all the special ams fonts etc.
\usepackage{graphicx}

\usepackage[authoryear,sort&compress]{natbib}
\usepackage{url}
\usepackage{enumerate}
\usepackage{colordvi,color,pspicture}
\usepackage{psfrag}

\usepackage{setspace}
\graphicspath{{Figures/}}

\newcommand{\Y}{\mathcal{Y}}
\newcommand{\X}{\mathcal{X}}
\newcommand{\U}{\mathcal{U}}
\newcommand{\I}{\mathbf{I}}
\newcommand{\C}{\mathcal{C}}
\newcommand{\mI}{\mathcal{I}}

\newcommand{\F}{\mathcal F}

\newcommand{\G}{\Gamma}
\newcommand{\N}{\mathcal{N}}
\newcommand{\V}{\mathcal{V}}
\newcommand{\A}{\mathcal{A}}
\newcommand{\R}{\mathbb{R}}
\newcommand{\mS}{\mathcal{S}}
\newcommand{\E}{\mathbf{E}}
\newcommand{\Var}{\mathbf{Var}}
\newcommand{\Cov}{\mathbf{Cov}}
\newcommand{\D}{\mathscr D}

\newcommand{\TY}{T(\Y_3)}

\newcommand{\y}{\mathsf{y}}

\theoremstyle{plain}
\newtheorem{theorem}{Theorem}[section]
\newtheorem{lemma}[theorem]{Lemma}

\theoremstyle{definition}

\theoremstyle{remark}

\newtheorem{remark}[theorem]{Remark}

\begin{document}

\title{Technical Report \# KU-EC-11-1:\\
Distribution of the Relative Density of Central Similarity Proximity Catch Digraphs
Based on One Dimensional Uniform Data}
\author{
Elvan Ceyhan\thanks{Address:
Department of Mathematics, Ko\c{c} University, 34450 Sar{\i}yer, Istanbul, Turkey.
e-mail: elceyhan@ku.edu.tr, tel:+90 (212) 338-1845, fax: +90 (212) 338-1559.
}
}

\date{\today}

\maketitle

{\bf short title:}
Relative Density of Central Similarity Proximity Catch Digraphs

\begin{abstract}
\noindent
We consider the distribution of a graph invariant
of central similarity proximity catch digraphs (PCDs)
based on one dimensional data.
The central similarity PCDs
are also a special type of parameterized
random digraph family
defined with two parameters,
a centrality parameter and an expansion parameter,
and for one dimensional data,
central similarity PCDs can also be viewed as a type of interval catch digraphs.
The graph invariant we consider is the relative density of central similarity PCDs.
We prove that relative density of central similarity PCDs is a $U$-statistic
and obtain the asymptotic normality
under mild regularity conditions using the central limit theory of $U$-statistics.
For one dimensional uniform data,
we provide the asymptotic distribution of the relative density of the central similarity PCDs
for the entire ranges of centrality and expansion parameters.
Consequently, we determine the optimal parameter values at which
the rate of convergence (to normality) is fastest.
We also provide the connection with class cover catch digraphs
and the extension of central similarity PCDs to higher dimensions.
\end{abstract}

%\vspace{0.1 in}



%\indent
%$^\star$
%This work was partially sponsored by TUBITAK Kariyer Project Grant 107T647.\\

\section{Introduction}
\label{sec:intro}
Proximity catch digraphs (PCDs) are introduced recently and have applications
in spatial data analysis and statistical pattern classification.
The PCDs are a special type of proximity graphs which were introduced by \cite{toussaint:1980}.
Furthermore, the PCDs are closely related to the class cover problem of \cite{cannon:2000}.
The PCDs are vertex-random digraphs in which each vertex corresponds to a data point,
and directed edges (i.e., arcs) are defined by some bivariate relation on the data
using the regions based on these data points.

\cite{priebe:2001} introduced the class cover catch digraphs (CCCDs) in $\R$
which is a special type of PCDs and gave the exact and
the asymptotic distribution of the domination number of the CCCDs
based on data from two classes, say $\X$ and $\Y$, with uniform distribution on a bounded interval in $\R$.
\cite{devinney:2002a}, \cite{marchette:2003}, \cite{priebe:2003b},
\cite{priebe:2003a}, and \cite{devinney:2006} applied the concept in higher dimensions and
demonstrated relatively good performance of CCCDs in classification.
\cite{ceyhan:CS-JSM-2003} introduced central similarity PCDs for two dimensional data in an unparameterized fashion;
the parameterized version of this PCD is later developed by \cite{ceyhan:arc-density-CS}
where the relative density of the PCD is calculated and used for testing bivariate spatial patterns in $\R^2$.
\cite{ceyhan:dom-num-NPE-SPL,ceyhan:dom-num-NPE-MASA}, \cite{ceyhan:dom-num-NPE-Spat2011}
applied the same concept (for a different PCD family called proportional-edge PCD)
in testing spatial point patterns in $\mathbb R^2$.
The distribution of the relative density of the proportional-edge PCDs for one dimensional uniform data
is provided in \cite{ceyhan:rel-dens-NPE-1D}.

In this article,
we consider central similarity PCDs
for one dimensional data.
We derive the asymptotic distribution of a graph invariant called
\emph{relative (arc) density} of central similarity PCDs.
Relative density is the ratio of number of arcs in a given digraph with $n$ vertices
to the total number of arcs possible (i.e., to the number of arcs in a complete symmetric digraph of order $n$).
We prove that, properly scaled,
the relative density of the central similarity PCDs is a $U$-statistic,
which yields the asymptotic normality by the general central limit theory of $U$-statistics.
Furthermore, we derive the explicit form of the asymptotic normal distribution
of the relative density of the PCDs
for uniform one dimensional $\X$ points
whose support being partitioned by class $\Y$ points.
We consider the entire ranges of the expansion and centrality parameters
and the asymptotic distribution is derived
as a function of these parameters
based on detailed calculations.
The relative density of central similarity PCDs is first investigated
for uniform data in one interval (in $\R$) and the analysis is generalized to
uniform data in multiple intervals.
These results can be used in applying the relative density
for testing spatial interaction between classes of one dimensional data.
Moreover, the behavior of the relative density in the one dimensional
case forms the foundation of our investigation and extension of the topic in higher dimensions.

We define the proximity catch digraphs and describe the central similarity PCDs in Section \ref{sec:vert-random-PCDs},
define their relative density and provide
preliminary results in Section \ref{sec:rel-dens-PCDs},
provide the distribution of the relative density for uniform data in one interval
in Section \ref{sec:relative-density-uniform}
and in multiple intervals in Section \ref{sec:dist-multiple-intervals},
provide extension to higher dimensions in Section \ref{sec:NCSt-higher-D}
and provide discussion and conclusions in Section \ref{sec:disc-conclusions}.
Shorter proofs are given in the main body of the article;
while longer proofs are deferred to the Appendix Sections.

\section{Vertex-Random Proximity Catch Digraphs}
\label{sec:vert-random-PCDs}
We first define vertex-random PCDs in a general setting.
Let $(\Omega,\mathcal M)$ be a measurable space and $\X_n=\{X_1,X_2,\ldots,X_n\}$ and
$\Y_m=\{Y_1,Y_2,\ldots,Y_m\}$ be two sets of $\Omega$-valued random variables
from classes $\X$ and $\Y$, respectively, with joint probability distribution $F_{X,Y}$
and marginals $F_X$ and $F_Y$, respectively.
A PCD is comprised by a set $\V$ of vertices and a set $\A$ of arcs.
For example, in the two class case,
with classes $\X$ and $\Y$,
we choose the $\X$ points to be the
vertices and put an arc from $X_i \in \X_n$ to $X_j\in \X_n$,
based on a binary relation which measures the relative
allocation of $X_i$ and $X_j$ with respect to $\Y$ points.
Notice that the randomness is only on the vertices,
hence the name \emph{vertex-random PCDs}.
Consider the map $N:\Omega \rightarrow \mathcal P(\Omega)$,
where $\mathcal P(\Omega)$ represents the power set of $\Omega$.
Then given $\Y_m \subseteq \Omega$,
the {\em proximity map}
$N(\cdot)$ associates with each point $x \in \Omega$
a {\em proximity region} $N(x) \subseteq \Omega$.
For $B \subseteq \Omega$, the $\G_1$-region is the image of the map
$\G_1(\cdot,N):\mathcal P(\Omega) \rightarrow \mathcal P(\Omega)$
that associates the region $\G_1(B,N):=\{z \in \Omega: B \subseteq  N(z)\}$
with the set $B$.
For a point $x \in \Omega$, we denote $\G_1(\{x\},N)$ as $\G_1(x,N)$.
Notice that while the proximity region is defined for one point,
a $\G_1$-region is defined for a point or set of points.
The {\em vertex-random PCD} has the vertex set $\V=\X_n$
and arc set $\A$ defined by $(X_i,X_j) \in \A$ if $X_j \in N(X_i)$.
Let arc probability be defined as
$p_a(i,j):=P((X_i,X_j) \in \A)$
for all $i \not=j$, $i,j =1,2,\ldots,n$.
Given $\Y_m=\{y_1,y_2,\ldots,y_m\}$, let $\X_n$ be a random sample from $F_X$.
Then $N(X_i)$ are also iid and the same holds for $\G_1(X_i,N)$.
Hence $p_a(i,j)=p_a$
for all $i \not=j$, $i,j =1,2,\ldots,n$ for such $\X_n$.

\subsection{Central Similarity PCDs for One Dimensional Data}
\label{sec:cent-sim-PCD}
\begin{figure} [h]
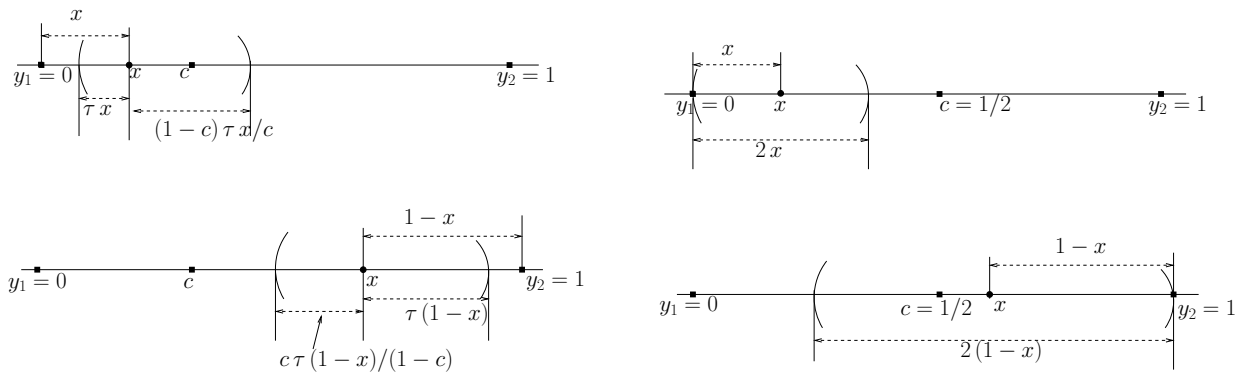

\centering
{\Huge
\scalebox{.35}{\input{NCSonedim.pstex_t}}
\hspace{1cm}
\scalebox{.35}{\input{CCCDonedim2.pstex_t}}
}
\caption
{Plotted in the left is an illustration of the construction of central similarity proximity region, $N(x,\tau,c)$
with $\tau \in (0,1)$, $\Y_2=\{y_1,y_2\}$ with $y_1=0$ and $y_2=1$ (hence $M_c=c$)
and $x \in (0,c)$ (top) and $x \in (c,1)$ (bottom);
and in the right is the
proximity region associated with CCCD, i.e., $N(x,\tau=1,c=1/2)$
for an $x \in (0,1/2)$ (top) and $x \in (1/2,1)$ (bottom).}
\label{fig:ProxMapDef1D}
\end{figure}

In the special case of central similarity PCDs for one dimensional data,
we have $\Omega=\R$.
Let $Y_{(i)}$ be the $i^{th}$ order statistic of $\Y_m$ for $i=1,2,\ldots,m$.
Assume $Y_{(i)}$ values are distinct (which happens with probability one for continuous distributions).
Then $Y_{(i)}$ values partition $\R$ into $(m+1)$ intervals.
Let
$$-\infty =: Y_{(0)}<Y_{(1)}< \ldots <Y_{(m)}< Y_{(m+1)}:=\infty.$$
We call intervals $(-\infty,Y_{(1)})$ and $\left(Y_{(m)},\infty\right)$ the \emph{end intervals},
and intervals $(Y_{(i-1)},Y_{(i)})$ for $i=2,\ldots,m$ the \emph{middle intervals}.
Then we define the central similarity PCD with the parameter $\tau > 0$
for two one dimensional data sets, $\X_n$ and $\Y_m$,
from classes $\X$ and $\Y$, respectively, as follows.
For $x \in \left( Y_{(i-1)},Y_{(i)} \right)$ with $i \in \{2,\ldots,m\}$
(i.e., for $x$ in a middle interval)
and $M_c \in \left( Y_{(i-1)},Y_{(i)} \right)$ such that $c \times 100$ \% of
$(Y_{(i)}-Y_{(i-1)})$ is to the left of $M_c$ (i.e., $M_c=Y_{(i-1)}+c\,(Y_{(i)}-Y_{(i-1)})$)
\begin{equation}
\label{eqn:NCSt-general-defn1}
N(x,\tau,c)=
\begin{cases}
\left(x-\tau\,\left(x-Y_{(i-1)}\right),x+\frac{\tau\,(1-c)\left(x-Y_{(i-1)}\right)}{c}\right) \bigcap \left( Y_{(i-1)},Y_{(i)} \right) & \text{if $x \in (Y_{(i-1)},M_c)$,}
\vspace{0.2cm}\\
\left(x-\frac{c\,\tau\,\left(Y_{(i)}-x\right)}{1-c},x+\tau\,\left(Y_{(i)}-x\right)\right)  \bigcap \left( Y_{(i-1)},Y_{(i)} \right)     & \text{if $x \in \left( M_c,Y_{(i)} \right)$.}
\end{cases}
\end{equation}
Observe that with $\tau \in (0,1)$,
we have
\begin{equation}
\label{eqn:NCSt-general-defn-t01}
N(x,\tau,c)=
\begin{cases}
\left(x-\tau\,\left(x-Y_{(i-1)}\right),x+\frac{\tau\,(1-c)\left(x-Y_{(i-1)}\right)}{c}\right) & \text{if $x \in (Y_{(i-1)},M_c)$,}
\vspace{0.2cm}\\
\left(x-\frac{c\,\tau\,\left(Y_{(i)}-x\right)}{1-c},x+\tau\,\left(Y_{(i)}-x\right)\right) & \text{if $x \in \left( M_c,Y_{(i)} \right)$,}
\end{cases}
\end{equation}
and with $\tau \ge 1$,
we have
\begin{equation}
\label{eqn:NCSt-general-defn-gtr1}
N(x,\tau,c)=
\begin{cases}
\left(Y_{(i-1)},x+\frac{\tau\,(1-c)\left(x-Y_{(i-1)}\right)}{c}\right) & \text{if $x \in \left(Y_{(i-1)},\frac{c\,Y_{(i)}+\tau\,(1-c)\,Y_{(i-1)}}{c+\tau\,(1-c)}\right)$,}
\vspace{0.2cm}\\
\left(Y_{(i-1)},Y_{(i)}\right) & \text{if $x \in \left(\frac{c\,Y_{(i)}+\tau\,(1-c)\,Y_{(i-1)}}{c+\tau\,(1-c)},\frac{(1-c)\,Y_{(i-1)}+c\,\tau\,Y_{(i)}}{1-c+c\,\tau}\right)$,}
\vspace{0.2cm}\\
\left(x-\frac{c\,\tau\,\left(Y_{(i)}-x\right)}{1-c},Y_{(i)}\right) & \text{if $x \in \left( \frac{(1-c)\,Y_{(i-1)}+c\,\tau\,Y_{(i)}}{1-c+c\,\tau},Y_{(i)} \right)$.}
\end{cases}
\end{equation}
For an illustration of $N(x,\tau,c)$ in the middle interval case,
see Figure \ref{fig:ProxMapDef1D} (left)
where $\Y_2=\{y_1,y_2\}$ with $y_1=0$ and $y_2=1$
(hence $M_c=c$).

Additionally,
for $x \in \left( Y_{(i-1)},Y_{(i)} \right)$ with $i \in \{1,m+1\}$
(i.e., for $x$ in an end interval),
the central similarity proximity region only has an expansion parameter,
but not a centrality parameter.
Hence we let $N_e(x,\tau)$ be the central similarity proximity region for an $x$ in an end interval.
Then
with $\tau \in (0,1)$,
we have
\begin{equation}
\label{eqn:NPEr-general-defn2-t01}
N_e(x,\tau)=
\begin{cases}
\left(x-\tau\,\left(Y_{(1)}-x\right),x+\tau\,\left(Y_{(1)}-x\right)\right)     & \text{if $x < Y_{(1)}$,}\\
\left(x-\tau\,\left(x-Y_{(m)}\right),x+\tau\,\left(x-Y_{(m)}\right)\right)  & \text{if $x > Y_{(m)}$}
\end{cases}
\end{equation}
and with $\tau \ge 1$,
we have
\begin{equation}
\label{eqn:NPEr-general-defn2-tgtr1}
N_e(x,\tau)=
\begin{cases}
\left(x-\tau\,\left(Y_{(1)}-x\right),Y_{(1)}\right)     & \text{if $x < Y_{(1)}$,}\\
\left(Y_{(m)},x+\tau\,\left(x-Y_{(m)}\right)\right)  & \text{if $x > Y_{(m)}$.}
\end{cases}
\end{equation}
%Notice that for $i \in \{1,m+1\}$,
%the region
%$N(x,\tau,c)$ is actually independent of $c$.
If $x \in \Y_m$, then we define $N(x,\tau,c)=\{x\}$ and $N_e(x,\tau)=\{x\}$ for all $\tau > 0$,
and if $x = M_c$, then in Equation \eqref{eqn:NCSt-general-defn1},
we arbitrarily assign $N(x,\tau,c)$ to be one of
$\left(x-\tau\,\left(x-Y_{(i-1)}\right),x+\frac{\tau\,(1-c)\left(x-Y_{(i-1)}\right)}{c}\right) \bigcap \left( Y_{(i-1)},Y_{(i)} \right) $ or
$\left(x-\frac{c\,\tau\,\left(Y_{(i)}-x\right)}{1-c},x+\tau\,\left(Y_{(i)}-x\right)\right)  \bigcap \left( Y_{(i-1)},Y_{(i)} \right)$.
For $X$ from a continuous distribution,
these special cases in the construction of central similarity proximity region
--- $X \in \Y_m$ and $X = M_c$ --- happen with probability zero.
Notice that $\tau > 0$ implies $x \in N(x,\tau,c)$ for all $x \in \left[ Y_{(i-1)},Y_{(i)} \right]$ with $i \in \{2,\ldots,m\}$
and $x \in N_e(x,\tau)$ for all $x \in \left[ Y_{(i-1)},Y_{(i)} \right]$ with $i \in \{1,m+1\}$.
Furthermore,
$\lim_{\tau \rightarrow \infty} N(x,\tau,c) = \left( Y_{(i-1)},Y_{(i)} \right)$
(and $\lim_{\tau \rightarrow \infty} N_e(x,\tau) = \left( Y_{(i-1)},Y_{(i)} \right)$)
for all $x \in \left( Y_{(i-1)},Y_{(i)} \right)$ with $i \in \{2,\ldots,m\}$ (and $i \in \{1,m+1\}$),
so we define $N(x,\infty,c) = \left( Y_{(i-1)},Y_{(i)} \right)$
(and $N_e(x,\infty) = \left( Y_{(i-1)},Y_{(i)} \right)$) for all such $x$.
%For $x \in \Y_m$, we define $N(x,\tau,c) = \{x\}$ for all $\tau > 0$.

The vertex-random central similarity PCD has the vertex set $\X_n$ and arc set $\A$ defined by
$(X_i,X_j) \in \A \iff X_j \in N(X_i,\tau,c)$ for $X_i,X_j$ in the middle intervals
and
$(X_i,X_j) \in \A \iff X_j \in N_e(X_i,\tau)$ for $X_i,X_j$ in the end intervals.
We denote such digraphs as $\D_{n,m}(\tau,c)$.
A $\D_{n,m}(\tau,c)$-digraph is a {\em pseudo digraph} according to some authors,
if loops are allowed (see, e.g., \cite{chartrand:1996}).
The $\D_{n,m}(\tau,c)$-digraphs are closely related to the {\em proximity graphs} of
\cite{jaromczyk:1992} and might be considered as a special case of
{\em covering sets} of \cite{tuza:1994}.
Our vertex-random proximity digraph is not a standard random graph (see, e.g., \cite{janson:2000}).
The randomness of a $\D_{n,m}(\tau,c)$-digraph lies in the fact that
the vertices are random with the joint distribution $F_{X,Y}$,
but arcs $(X_i,X_j)$ are
deterministic functions of the random variable $X_j$ and the random set $N(X_i,\tau,c)$ in the middle intervals
and
the random set $N_e(X_i,\tau)$ in the end intervals.
In $\R$, the vertex-random PCD is a special case of
{\em interval catch digraphs} (see, e.g., \cite{sen:1989} and \cite{prisner:1994}).
Furthermore, when $\tau=1$ and $c=1/2$ (i.e., $M_c=\left( Y_{(i-1)}+Y_{(i)} \right)/2$)
we have $N(x,1,1/2)=B(x,r(x))$ for an $x$ in a middle interval
and
$N_e(x,1)=B(x,r(x))$ for an $x$ in an end interval
where $r(x)=d(x,\Y_m)=\min_{y \in \Y_m}d(x,y)$
and the corresponding PCD is the CCCD of \cite{priebe:2001}.
See also Figure \ref{fig:ProxMapDef1D} (right).

\section{Relative Density of Vertex-Random PCDs}
\label{sec:rel-dens-PCDs}
Let $D_n=(\V,\A)$ be a digraph with vertex set $\V=\{v_1,v_2,\ldots,v_n\}$
and arc set $\A$
and let $|\cdot|$ stand for the set cardinality function.
The relative density of the digraph $D_n$ which is of order $|\V| = n \ge 2$,
denoted $\rho(D_n)$, is defined as (\cite{janson:2000})
$$\rho(D_n) = \frac{|\A|}{n(n-1)}.$$
Thus $\rho(D_n)$ represents the ratio of the number of arcs
in the digraph $D_n$ to the number of arcs in the complete symmetric digraph of order $n$,
which is $n(n-1)$.
For $n \le 1$, we set $\rho(D_n)=0$,
since there are no arcs.
If $D_n$ is a random digraph in which arcs %and/or vertices
result from a random process,
then the \emph{arc probability} between vertices $v_i,v_j$ is
$p_a(i,j)=P((v_i,v_j) \in \A)$ for all $i \not=j$, $i,j =1,2,\ldots,n$.

Given $\Y_m=\{y_1,y_2,\ldots,y_m\}$,
let $\X_n$ be a random sample from $F_X$
and $D_n$ be the PCD based on proximity region $N(\cdot)$ with vertices $\X_n$
and the arc set $\A$ is defined as $(X_i,X_j) \in \A$ if $X_j \in N(X_i)$.
Let $h_{ij} := (g_{ij}+g_{ji})/2$ where $g_{ij}=I((X_i,X_j) \in \A)+\I(X_j \in N(X_i))$.
Then we can rewrite the relative density as follows:
$$
\rho(D_n)=
\frac{2}{n(n-1)}
\sum\hspace*{-0.1 in}\sum_{i < j \hspace*{0.25 in}}   \hspace*{-0.1 in} h_{ij}.
$$
Although the digraph is asymmetric,
$h_{ij}$ is defined as the average number of arcs between $X_i$ and $X_j$
in order to produce a symmetric kernel with finite variance (\cite{lehmann:1988}).
The relative density
$\rho(D_n)$ is a random variable that depends on $n$, $F$, and $N(\cdot)$ (i.e., $\Y_m$).
But $\E\left[ \rho(D_n) \right]=\E\left[ h_{12} \right]=p_a$
only depends on $F$ and $N(\cdot)$.
Furthermore,
\begin{equation}
\label{eqn:Var[rho-and-D]}
0 \le \Var\left[ \rho(D_n) \right]=
\frac{4}{n^2\,(n-1)^2}\Var\left[\sum\hspace*{-0.1 in}\sum_{i < j \hspace*{0.25 in}}   \hspace*{-0.1 in} h_{ij}\right]
=\frac{2}{n\,(n-1)}\Var\left[ h_{12} \right]+
\frac{4(n-2)}{n\,(n-1)} \, \Cov\left[ h_{12},h_{13} \right] \le 1/4.
\end{equation}
Hence $\rho(D_n)$, is a one-sample $U$-statistic of degree 2
and is an unbiased estimator of arc probability $p_a$.
If, additionally, $\nu=\Cov\left[ h_{ij},h_{ik} \right]>0$ for all $i \not=j\not=k$,
$i,j,k \in \{1,2,\ldots,n\}$,
then a CLT for $U$-statistics (\cite{lehmann:1988}) yields
$\sqrt{n}\,\left[\rho(D_n)-p_a\right]\stackrel{\mathcal L}{\longrightarrow} \N(0,4\,\nu)$
as $n \rightarrow \infty$,
where $\stackrel{\mathcal L}{\longrightarrow}$ stands for convergence in law
and
$\N(\mu,\sigma^2)$ stands for the normal distribution with mean $\mu$ and variance $\sigma^2$.

In Equation \eqref{eqn:Var[rho-and-D]}, we have
\begin{multline*}
\Var[h_{ij}]=\Var\left[ h_{12} \right]=
\E[(h_{ij})^2]-(\E[h_{ij}])^2=
\E[(g_{ij}+g_{ji})^2/4]-p_a^2=\\
(\E[g_{ij}]+2\,\E[g_{ij}]\E[g_{ji}]+\E[g_{ji}])/4-p_a^2=
(p_a+2\,p_a+p_a)/4-p_a^2=
(p_a-p_a^2)/2=
p_a\left(1-p_a\right)/2
\end{multline*}
and the covariance is
$$\Cov\left[ h_{12},h_{13} \right]=
\E\left[ h_{12}h_{13} \right]-\E\left[ h_{12} \right]\E\left[ h_{13} \right]=
\E\left[ h_{12}h_{13} \right]-p^2_a,$$
with
\begin{eqnarray*}
4\,\E\left[ h_{12}h_{13} \right]& = & \E[(g_{12}+g_{21})(g_{13}+g_{31})]=\E[g_{12}g_{13}+g_{12}g_{31}+g_{21}g_{13}+g_{21}g_{31}]\\
& = & \E[\I(X_2 \in N(X_1)\I( X_3 \in N(X_1))+\I(X_2 \in N(X_1)\I( X_1 \in N(X_3))+\\
&  & \I(X_1 \in N(X_2)\I( X_3 \in N(X_1))]+\I(X_1 \in N(X_2)\I( X_1 \in N(X_3))]\\
& = & \E[\I(\{X_2,X_3\} \subset N(X_1))+\I(X_2 \in N(X_1)\I(X_3 \in \G_1(X_3,N))+\\
& & \I(X_2 \in \G_1(X_1)\I( X_3 \in N(X_1))]+\I(X_2 \in \G_1(X_1,N)\I( X_3 \in \G_1(X_1,N))]\\
& = & P(\{X_2,X_3\} \subset N(X_1))+2\,P(X_2 \in N(X_1),X_3 \in \G_1(X_1,N))+P(\{X_2,X_3\} \subset \G_1(X_1,N)).
\end{eqnarray*}

Then $\nu=\Cov(h_{ij},h_{ik})=\E[h_{ij}h_{ik}]-\E[h_{ij}]\E[h_{ik}]=
\E[h_{ij}h_{ik}]-p_a^2=\E[h_{12}h_{13}]-p_a^2>0$ iff
$$P(\{X_2,X_3\} \subset N(X_1))+2\,P(X_2 \in N(X_1),X_3 \in \G_1(X_1,N))+P(\{X_2,X_3\} \subset \G_1(X_1,N))> 4 p_a^2.$$

Notice also that
\begin{multline*}
\E[|h_{ij}|^3]=
\E[(g_{ij}+g_{ji})^3/8]=
\E[g^3_{ij}+3\,g^2_{ij}g_{ji}+3\,g_{ij}g^2_{ji}+g^3_{ji}]/8=
\E[g_{ij}+3\,g_{ij}g_{ji}+3\,g_{ij}g_{ji}+g_{ji}]/8=\\
(2\,\E[g_{ij}]+6\,\E[g_{ij}]\E[g_{ji}])/8=
(p_a+3\,p^2_a)/4 < \infty.
\end{multline*}
Then for $\nu >0$,
the sharpest rate of convergence in the asymptotic normality of
$\rho(D_n)$ is
\begin{equation}
\label{eqn:rate-of-convergence}
\sup_{t\in \R} \left| P \left( \frac{\sqrt{n}(\rho(D_n)-p_a)}{\sqrt{4\,\nu}}\le t \right)-\Phi(t) \right|
\le 8\,K\,p_a\, (4\,\nu)^{-3/2}\, n^{-1/2}=
K\,\frac{p_a}{\sqrt{n\,\nu^3}}
\end{equation}
where $K$ is a constant and $\Phi(t)$ is the distribution function for
the standard normal distribution (\cite{callaert:1978}).

In general a random digraph, just like a random graph,
can be obtained by starting with a set of $n$ vertices and adding arcs between them at random.
We can consider the digraph counterpart of the Erd\H{o}s–--R\'{e}nyi model for random graphs,
denoted $D(n,p)$, in which every possible arc occurs independently with probability $p$ (\cite{erdos:1959}).
Notice that for the random digraph $D(n,p)$, the relative density of $D(n,p)$ is a $U$-statistic;
however, the asymptotic distribution of its relative density is degenerate
(with $\rho(D(n,p))\stackrel{\mathcal L}{\longrightarrow} p$, as $n \rightarrow \infty$)
since the covariance term is zero due to the independence between the arcs.

%\subsection{Relative Density of Random $\D_{n,m}(\tau,c)$-digraphs}
%\label{sec:rel-dens-digraph}
Let $\F(\mathbb R):=\{F_{X,Y} \text{ on } \mathbb R \text { with } P(X=Y)=0
\text{ and the marginals, } \text{$F_X$ and $F_Y$, are non-atomic}\}$.
In this article,
we consider $\D_{n,m}(\tau,c)$-digraphs for which
$\X_n$ and $\Y_m$ are random samples from $F_X$ and $F_Y$, respectively,
and the joint distribution of $X,Y$ is $F_{X,Y} \in \F(\R)$.
Then the order statistics of $\X_n$ and $\Y_m$ are distinct with probability one.
We call such digraphs as \emph{$\F(\R)$-random $\D_{n,m}(\tau,c)$-digraphs}
and focus on the random variable $\rho(\D_{n,m}(\tau,c))$.
%To make the dependence on sample sizes $n$ and $m$ and parameters $\tau$ and $c$, explicit,
For notational brevity,
we use $\rho_{n,m}(\tau,c)$ instead of $\rho(\D_{n,m}(\tau,c))$.
It is trivial to see that $0 \le \rho_{n,m}(\tau,c) \le 1$,
and $\rho_{n,m}(\tau,c) >0 $ for nontrivial digraphs.

\subsection{The Distribution of the Relative Density of $\F(\R)$-random $\D_{n,m}(\tau,c)$-digraphs}
\label{sec:relative-density-Dnm}
Let $\mI_i:=\left(Y_{(i-1)},Y_{(i)}\right)$, $\X_{[i]}:=\X_n \cap \mI_i$,
and $\Y_{[i]}:=\{Y_{(i-1)},Y_{(i)}\}$ for $i=1,2,\ldots,(m+1)$.
Let $D_{[i]}(\tau,c)$ be the component of the random $\D_{n,m}(\tau,c)$-digraph
induced by the pair $\X_{[i]}$ and $\Y_{[i]}$.
Then we have a disconnected digraph with subdigraphs $D_{[i]}(\tau,c)$
for $i=1,2,\ldots,(m+1)$
each of which might be null or itself disconnected.
Let $\A_{[i]}$ be the arc set of $D_{[i]}(\tau,c)$,
and $\rho_{{}_{[i]}}(\tau,c)$ denote the relative density
of $D_{[i]}(\tau,c)$;
$n_i:=\left|\X_{[i]}\right|$, and $F_i$ be the density $F_X$ restricted to $\mI_i$
for $i \in \{1,2,\ldots,m+1\}$.
Furthermore,
let $M^{[i]}_c \in \mI_i$ be the point so that it divides
the interval $\mI_i$ in ratios $c$ and $1-c$
(i.e., length of the subinterval to the left of $M^{[i]}_c$ is
$c \times 100$ \% of the length of $\mI_i$) for $i \in \{2,\ldots,m\}$.
Notice that for $i \in \{2,\ldots,m\}$ (i.e., middle intervals),
$D_{[i]}(\tau,c)$ is based on the proximity region $N(x,\tau,c)$
and
for $i \in \{1,m+1\}$ (i.e., end intervals),
$D_{[i]}(\tau,c)$ is based on the proximity region $N_e(x,\tau)$.
Since we have at most $m+1$ subdigraphs that are disconnected,
it follows that we have at most $n_{_T}:=\sum_{i=1}^{m+1} n_i (n_i-1)$ arcs
in the digraph $\D_{n,m}(\tau,c)$.
Then we define the relative density
for the entire digraph as
\begin{equation}
\label{eqn:rho-nm-tau,c}
\rho_{n,m}(\tau,c):=\frac{\left|\A\right|}{n_{_T}}
=\frac{\sum_{i=1}^{m+1} |\A_{[i]}|}{n_{_T}}
=\frac{1}{n_{_T}}\sum_{i=1}^{m+1}(n_i (n_i-1))\rho_{{}_{[i]}}(\tau,c).
\end{equation}
Since $\frac{n_i\,(n_i-1)}{n_{_T}} \ge 0$ for each $i$
and $\displaystyle \sum_{i=1}^{m+1}\frac{n_i\,(n_i-1)}{n_{_T}}=1$,
it follows that
$\rho_{n,m}(\tau,c)$ is a mixture of the $\rho_{{}_{[i]}}(\tau,c)$.
We study the simpler random variable $\rho_{{}_{[i]}}(\tau,c)$ first.
In the remaining of this section,
the almost sure (a.s.) results follow from the fact that
the marginal distributions $F_X$ and $F_Y$ are non-atomic.

\begin{lemma}
\label{lem:end-intervals}
Let $D_{[i]}(\tau,c)$ be the digraph induced by $\X$ points in the end intervals
(i.e., $i \in\{ 1,(m+1)\}$) and $\rho_{{}_{[i]}}(\tau,c)$ be the corresponding relative density.
For $\tau > 0$,
if $n_i \le 1$, then $\rho_{{}_{[i]}}(\tau,c)=0$.
For $\tau \ge 1$,
if $n_i >1$, then $\rho_{{}_{[i]}}(\tau,c)\ge 1/2$ a.s.
\end{lemma}

\noindent
{\bf Proof:}
Let $i=m+1$ (i.e., consider the right end interval).
For all $\tau > 0$,
if $n_{m+1} \le 1$, then by definition $\rho_{{}_{[m+1]}}(\tau,c)=0$.
So, we assume $n_{m+1}>1$.
Let $\X_{[m+1]}=\{Z_1,Z_2,\ldots,Z_{n_{m+1}}\}$
and $Z_{(j)}$ be the corresponding order statistics.
Then for $\tau \ge 1$,
there is an arc from $Z_{(j)}$ to each $Z_{(k)}$ for $k<j$,
with $j,k \in\{1,2,\ldots,n_{m+1}\}$
(and possibly to some other $Z_l$),
since $N_e\left(Z_{(j)},\tau\right)=(Y_{(m)},Z_{(j)}+\tau\,(Z_{(j)}-Y_{(m)}))$
and so $Z_{(k)} \in N_e\left(Z_{(j)},\tau\right)$.
So, there are at least $0+1+2+\ldots+n_{m+1}-1=n_{m+1}(n_{m+1}-1)/2$ arcs in $D_{[m+1]}(\tau,c)$.
Then $\rho_{{}_{[i]}}(\tau,c) \ge (n_{m+1}(n_{m+1}-1)/2)/(n_{m+1}(n_{m+1}-1))=1/2$.
By symmetry, the same results hold for $i=1$.
$\blacksquare$

Using Lemma \ref{lem:end-intervals},
we obtain the following lower bound for $\rho_{n,m}(\tau,c)$ for $\tau \ge 1$.

\begin{theorem}
\label{thm:rho-Dnm-r-M}
Let $D_{n,m}(\tau,c)$ be an $\F(\R)$-random $\D_{n,m}(\tau,c)$-digraph with $n>0,\,m>0$
and $k_1$ and $k_2$ be two natural numbers defined as
$k_1:=\sum_{i=2}^m (n_{i,1}(n_{i,1}-1)/2+n_{i,2}(n_{i,2}-1)/2)$
and
$k_2:=\sum_{i \in \{1,m+1\}} n_i(n_i-1)/2$,
where
$n_{i,1}:=\left|\X_n \cap \left( Y_{(i-1)},M^{[i]}_c \right)\right|$
and
$n_{i,2}:=\left|\X_n \cap \left(M^{[i]}_c,Y_{(i)}\right)\right|$.
Then for $\tau \ge 1$,
we have $(k_1+k_2)/n_{_T} \le \rho_{n,m}(\tau,c)\le 1$ a.s.
\end{theorem}

\noindent
{\bf Proof:}
For $i \in\{ 1,(m+1)\}$,
we have $k_2$ as in Lemma \ref{lem:end-intervals}.
Let $i \in \{2,3,\ldots,m\}$
and
$\X_{i,1}:=\X_{[i]} \cap \left( Y_{(i-1)},M^{[i]}_c \right)=\{U_1,U_2,\ldots,U_{n_{i,1}}\}$,
and
$\X_{i,2}:=\X_{[i]} \cap \left(M^{[i]}_c,Y_{(i)}\right)=\{V_1,V_2,\ldots,V_{n_{i,2}}\}$.
Furthermore,
let $U_{(j)}$ and $V_{(k)}$ be the corresponding order statistics.
For $\tau \ge 1$,
there is an arc from $U_{(j)}$ to $U_{(k)}$ for $k < j$, $j,k \in\{1,2,\ldots,n_{i,1}\}$
and possibly to some other $U_l$,
and similarly
there is an arc from $V_{(j)}$ to $V_{(k)}$ for $k > j$, $j,k \in\{1,2,\ldots,n_{i,2}\}$
and possibly to some other $V_l$.
Thus there are at least
$\frac{n_{i,1}(n_{i,1}-1)}{2}+\frac{n_{i,2}(n_{i,2}-1)}{2}$ arcs in $D_{[i]}(\tau,c)$.
Hence $\rho_{n,m}(\tau,c) \ge (k_1+k_2)/n_{_T}$.
$\blacksquare$

\begin{theorem}
\label{thm:rho-for-r=infty}
For $i=1,2,3,\ldots,m+1$,
$\tau=\infty$, and $n_i > 0$,
we have
$\rho_{{}_{[i]}}(\tau=\infty,c) =\I(n_i > 1)$
and
$\rho_{n,m}(\tau=\infty,c)=1$ a.s.
\end{theorem}

\noindent
{\bf Proof:}
For $\tau=\infty$,
if $n_i \le 1$,
then $\rho_{{}_{[i]}}(\tau=\infty,c) =0$.
So we assume $n_i >1 $ and let $i=m+1$.
Then $N_e(x,\infty)=\left(Y_{(m)},\infty\right)$ for all $x \in \left(Y_{(m)},\infty\right)$.
Hence $D_{[m+1]}(\infty,c)$ is a complete symmetric digraph of order $n_{m+1}$,
which implies $\rho_{{}_{[m+1]}}(\tau=\infty,c)=1$.
By symmetry, the same holds for $i=1$.
For $i \in \{2,3,\ldots,m\}$ and $n_i >1 $,
we have $N(x,\infty,c)=\mI_i$ for all $x \in \mI_i$,
hence $D_{[i]}(\infty,c)$ is a complete symmetric digraph of order $n_{i}$,
which implies $\rho_{{}_{[i]}}(\infty,c)=1$.
Then
$\rho_{n,m}(\infty,c)=\sum \frac{n_i(n_i-1)\rho_{{}_{[i]}}(\infty,c)}{n_{_T}}=1$,
since when $n_i \le 1$, $n_i$ has no contribution to $n_{_T}$,
and when $n_i > 1$, we have $\rho_{{}_{[i]}}(\infty,c)=1$.
$\blacksquare$

\section{The Distribution of the Relative Density of Central Similarity PCDs
for Uniform Data}
\label{sec:relative-density-uniform}
Let $-\infty<\delta_1<\delta_2<\infty$,
$\Y_m$ be a random sample from non-atomic $F_Y$ with support $\mS(F_Y) \subseteq (\delta_1,\delta_2)$,
and $\X_n =\{X_1,X_2,\ldots,X_n\}$ be a random sample from $F_X=\U(\delta_1,\delta_2)$,
the uniform distribution on $(\delta_1,\delta_2)$.
So we have $F_{X,Y} \in \F(\R)$.
Assuming we have the realization of $\Y_m$ as
$\Y_m=\{y_1,y_2,\ldots,y_m\}=\{y_{(1)},y_{(2)},\ldots,y_{(m)}\}$ with $\delta_1<y_{(1)}<y_{(2)}<\ldots<y_{(m)}<\delta_2$,
we let $y_{(0)}:=\delta_1$ and $y_{(m+1)}:=\delta_2$.
Then it follows that
the distribution of $X_i$ restricted to $\mI_i$ %, denoted as $F_X|_{\left(\delta_1,y_{(1)}\right)}$,
is $F_X|_{\mI_i}=\U(\mI_i)$.
%Similarly we have
%$F_X|_{(y_1,y_2)}=\U(y_1,y_2)$ and $F_X|_{\left(y_{(m)},\delta_2\right)}=\U\left(y_{(m)},\delta_2\right)$.
%Then $\mI_1=\left(\delta_1,y_{(1)}\right)$, $\mI_2=(y_1,y_2)$, and $\mI_3=\left(y_{(m)},\delta_2\right)$.
%In each of these three intervals, we have a $\D_{n,2}(\tau,c)$-digraph.
We call such digraphs as \emph{$\U(\delta_1,\delta_2)$-random $\D_{n,m}(\tau,c)$-digraphs}
and provide the distribution of their relative density
for the whole range of $\tau$ and $c$.
%Let $\rho_n(\tau,c)$ be the relative density
%of the PCD based on $N(\cdot,\tau,c)$ with $\X_n$.
We first present a ``scale invariance" result for central similarity PCDs.
This invariance property will simplify the notation in
our subsequent analysis by allowing us to consider the special case
of the unit interval $(0,1)$.

\begin{theorem}
\label{thm:scale-inv-NCSt}
(Scale Invariance Property)
Suppose $\X_n$ is a set of iid random variables from $\U(\delta_1,\delta_2)$
where $\delta_1<\delta_2$ and $\Y_m$ is set of $m$ distinct $\Y$ points in $(\delta_1,\delta_2)$.
Then for any $\tau > 0$, the distribution of $\rho_{{}_{[i]}}(\tau,c)$ is
independent of $\Y_{[i]}$ (and hence of the restricted support interval $\mI_i$) for all $i \in \{1,2,\ldots,m+1\}$.
\end{theorem}

\noindent
{\bf Proof:}
Let $\delta_1<\delta_2$ and $\Y_m$ be as in the hypothesis.
Any $\U(\delta_1,\delta_2)$ random variable can be transformed into a $\U(0,1)$
random variable by $\phi(x)=(x-\delta_1)/(\delta_2-\delta_1)$,
which maps intervals $(t_1,t_2) \subseteq (\delta_1,\delta_2)$ to
intervals $\bigl( \phi(t_1),\phi(t_2) \bigr) \subseteq (0,1)$.
That is, if $X \sim \U(\delta_1,\delta_2)$, then we have $\phi(X) \sim \U(0,1)$
and
$P(X \in (t_1,t_2))=P(\phi(X) \in \bigl( \phi(t_1),\phi(t_2) \bigr)$
for all $(t_1,t_2) \subseteq (\delta_1,\delta_2)$.
The distribution of $\rho_{{}_{[i]}}(\tau,c)$ is obtained by calculating such probabilities.
So, without loss of generality, we can assume $\X_{[i]}$
is a set of iid random variables from the $\U(0,1)$ distribution.
That is, the distribution of $\rho_{{}_{[i]}}(\tau,c)$ does not depend on $\Y_{[i]}$
and hence does not depend on the restricted support interval $\mI_i$.
$\blacksquare$

Note that scale invariance of $\rho_{{}_{[i]}}(\tau=\infty,c)$ follows trivially
for all $\X_n$ from any $F_X$ with support in $(\delta_1,\delta_2)$ with $\delta_1<\delta_2$,
since for $\tau=\infty$, we have $\rho_{{}_{[i]}}(\tau=\infty,c)=1$ a.s. for non-atomic $F_X$.

Based on Theorem \ref{thm:scale-inv-NCSt},
we may assume each $\mI_i$
as the unit interval $(0,1)$ for uniform data.
Then the central similarity proximity region for $x \in (0,1)$
with parameters $c \in (0,1)$ and $\tau > 0$
have the following forms.
If $x\in \mI_i$ for $i \in \{2,\ldots,m\}$ (i.e., in the middle intervals),
when transformed under $\phi(\cdot)$ to $(0,1)$, we have
\begin{equation}
\label{eqn:NCSt-(0,1)-defn1}
N(x,\tau,c)=
\begin{cases}
(x\,(1-\tau), x\,(c+(1-c)\,\tau)/c) \cap (0,1) & \text{if $x \in (0,c)$,}\\
(x-c\,\tau\,(1-x)/(1-c), x+(1-x)\,\tau)  \cap (0,1)     & \text{if $x \in (c,1)$.}
\end{cases}
\end{equation}
In particular,
for $\tau \in (0,1)$,
we have
\begin{equation}
\label{eqn:NCSt-(0,1)-defn1-t01}
N(x,\tau,c)=
\begin{cases}
(x\,(1-\tau), x\,(c+(1-c)\,\tau)/c) & \text{if $x \in (0,c)$,}\\
(x-c\,\tau\,(1-x)/(1-c), x+(1-x)\,\tau)     & \text{if $x \in (c,1)$}
\end{cases}
\end{equation}
and
for $\tau \ge 1$,
we have
\begin{equation}
\label{eqn:NCSt-(0,1)-defn1-tgtr1}
N(x,\tau,c)=
\begin{cases}
(0, x\,(c+(1-c)\,\tau)/c) & \text{if $x \in \left(0,\frac{c}{c+(1-c)\,\tau}\right)$,} \vspace{.1cm}\\
(0,1) & \text{if $x \in \left(\frac{c}{c+(1-c)\,\tau},\frac{c\,\tau}{1-c+c\,\tau}\right)$,}\vspace{.1cm}\\
(x-c\,\tau\,(1-x)/(1-c), 1)     & \text{if $x \in \left(\frac{c\,\tau}{1-c+c\,\tau},1\right)$.}
\end{cases}
\end{equation}
and $N(x=c,\tau,c)$ is arbitrarily taken to be one of
$(x\,(1-\tau), x\,(c+(1-c)\,\tau)/c) \cap (0,1)$ or $(x-c\,\tau\,(1-x)/(1-c), x+(1-x)\,\tau)  \cap (0,1)$.
This special case of ``$X=c$" happens with probability zero for uniform $X$.

If $x\in \mI_1$ (i.e., in the left end interval),
when transformed under $\phi(\cdot)$ to $(0,1)$, we have
$N_e(x,\tau)=(\max(0,x-\tau\,(1-x)),\min(1,x+\tau\,(1-x))$;
and
if $x\in \mI_{m+1}$ (i.e., in the right end interval),
when transformed under $\phi(\cdot)$ to $(0,1)$, we have
$N_e(x,\tau)=(\max(0,x\,(1-\tau)),\min(1,x\,(1+\tau)))$.

Notice that each subdigraph $D_{[i]}(\tau,c)$ is itself a
$\U(\mI_i)$-random $\D_{n,2}(\tau,c)$-digraph.
The distribution of the relative density of $D_{[i]}(\tau,c)$
is given in the following result.
\begin{theorem}
\label{thm:U(y1,y2-randomDn2(tau,c)}
Let $\rho_{{}_{[i]}}(\tau,c)$ be the relative density of
subdigraph $D_{[i]}(\tau,c)$ of the central similarity PCD based on uniform data in $(\delta_1,\delta_2)$
where $\delta_1<\delta_2$ and $\Y_m$ be a set of $m$ distinct $\Y$ points in $(\delta_1,\delta_2)$.
Then
for $\tau \in (0,\infty)$, as $n_i \rightarrow \infty$,
we have
\begin{itemize}
\item[(i)]
for $i \in \{2,\ldots,m\}$,
$\sqrt{n_i}\,\left[\rho_{{}_{[i]}}(\tau,c)-\mu(\tau,c)\right]\stackrel{\mathcal L}{\longrightarrow} \N(0,4\,\nu(\tau,c))$,
where $\mu(\tau,c)=\E\left[\rho_{{}_{[i]}}(\tau,c)\right]$ is the arc probability
and
$\nu(\tau,c)=\Cov[h_{12},h_{12}]$ in the middle intervals, and
\item[(ii)]
for $i \in \{1,m+1\}$,
$\sqrt{n_i}\,\left[\rho_{{}_{[i]}}(\tau,c)-\mu_e(\tau)\right]\stackrel{\mathcal L}{\longrightarrow} \N(0,4\,\nu_e(\tau))$,
where $\mu_e(\tau)=\E\left[\rho_{{}_{[i]}}(\tau,c)\right]$ is the arc probability
and
$\nu_e(\tau)=\Cov[h_{12},h_{12}]$ in the end intervals.
\end{itemize}
\end{theorem}

\noindent \textbf{Proof:}
(i)
Let $i \in \{2,\ldots,m\}$ (i.e., $\mI_i$ be a middle interval).
By the scale invariance for uniform data (see Theorem \ref{thm:scale-inv-NCSt}),
a middle interval can be assumed to be the unit interval $(0,1)$.
The mean of the asymptotic distribution of $\rho_{{}_{[i]}}(\tau,c)$ is computed as follows.
$$\E[\rho_{{}_{[i]}}(\tau,c)]=\E[h_{12}]=P(X_2 \in N(X_1,\tau,c))=\mu(\tau,c)$$
which is the arc probability.
And the asymptotic variance of $\rho_{{}_{[i]}}(\tau,c)$ is
$\Cov[h_{12},h_{13}]=4\,\nu(\tau,c)$.
For $\tau \in (0,\infty)$,
since $2h_{12} = \I(X_2 \in N(X_1,\tau,c))+\I(X_1 \in N(X_2,\tau,c))$
is the number of arcs between $X_1$ and $X_2$ in the PCD,
$h_{12}$ tends to be high if the proximity region $N(X_1,\tau,c)$ is large.
In such a case, $h_{13}$ tends to be high also.
That is, $h_{12}$ and $h_{13}$
tend to be high and low together.
So, for $\tau \in (0,\infty)$,
we have $\nu(\tau,c)>0$.
Hence asymptotic normality follows.
%See also Figures \ref{fig:asy-var-tau=1-c}, \ref{fig:asymean-var-tau-c=1/2},
%\ref{fig:asymean-var-tau-c}.

(ii)
In an end interval,
the mean of the asymptotic distribution of $\rho_{{}_{[i]}}(\tau,c)$ is
$$\E[\rho_{{}_{[i]}}(\tau,c)]=\E[h_{12}]=P(X_2 \in N_e(X_1,\tau))=\mu_e(\tau)$$
the asymptotic variance of $\rho_{{}_{[i]}}(\tau,c)$ is
$\Cov[h_{12},h_{13}]=4\,\nu_e(\tau)$.
For $\tau \in (0,\infty)$, as in (i),
we have $\nu_e(\tau)>0$.
Hence asymptotic normality follows.
%, and \ref{fig:asymean-var-tau-end} (right)
%and Appendix 1.
$\blacksquare$

Let
$P_{2N}:=P(\{X_2,X_3\} \subset N(X_1,\tau,c))$,
$P_{NG}:=P(X_2 \in N(X_1,\tau,c),X_3 \in \G_1(X_1,\tau,c))$,
and
$P_{2G}:=P(\{X_2,X_3\} \subset \G_1(X_1,\tau,c))$.
Then
$$
\Cov[h_{12},h_{13}]=\E[h_{12}h_{13}]-\E[h_{12}]\E[h_{13}]=\E[h_{12}h_{13}]-\mu(\tau,c)^2=(P_{2N}+2\,P_{NG}+P_{2G})/4-\mu(\tau,c)^2,
$$
since
\begin{multline*}
4\,\E[h_{12}h_{13}]=
P(\{X_2,X_3\} \subset N(X_1,\tau,c))+
2\,P(X_2 \in N(X_1,\tau,c),X_3 \in \G_1(X_1,\tau,c))+\\
P(\{X_2,X_3\} \subset \G_1(X_1,\tau,c))
=P_{2N}+2\,P_{NG}+P_{2G}.
\end{multline*}

Similarly, let
$P_{2N,e}:=P(\{X_2,X_3\} \subset N_e(X_1,\tau))$,
$P_{NG,e}:=P(X_2 \in N_e(X_1,\tau),X_3 \in \G_{1,e}(X_1,\tau))$,
and
$P_{2G,e}:=P(\{X_2,X_3\} \subset \G_{1,e}(X_1,\tau))$.
Then
$$
\Cov[h_{12},h_{13}]=(P_{2N,e}+2\,P_{NG,e}+P_{2G,e})/4-\mu_e(\tau)^2.
$$

For $\tau=\infty$,
we have
$N(x,\infty,c) =\mI_i$ for all $x \in \mI_i$ with $i \in \{2,\ldots,m\}$
and
$N_e(x,\infty) =\mI_i$ for all $x \in \mI_i$ with $i \in \{1,m+1\}$.
%Let $i \in \{1,2,3\}$.
Then for $i \in \{2,\ldots,m\}$
$$\E\left[ \rho_{{}_{[i]}}(\infty,c) \right] =
\E\left[ h_{12} \right]=
\mu(\infty,c)=
P(X_2 \in N(X_1,\infty,c)= P(X_2 \in \mI_i)=1.$$
On the other hand,
$4\,\E\left[ h_{12}h_{13} \right] =
P(\{X_2,X_3\} \subset N(X_1,\infty,c))+
2\,P(X_2 \in N(X_1,\infty,c),X_3 \in \G_1(X_1,\infty,c))+
P(\{X_2,X_3\} \subset \G_1(X_1,\infty,c))=
(1+2+1)$.
Hence $\E\left[ h_{12}h_{13} \right]=1$
and so
$\nu(\infty,c)=0$.
Similarly,
for $i \in \{1,m+1\}$,
we have
$\mu_e(\infty)=1$
and
$\nu_e(\infty)=0$.
Therefore, the CLT result does not hold for $\tau = \infty$.
Furthermore, $\rho_{{}_{[i]}}(\tau=\infty,c)=1$ a.s.

By Theorem \ref{thm:U(y1,y2-randomDn2(tau,c)},
we have $\nu(\tau,c) > 0$ (and $\nu_e(\tau) > 0$) iff
$P_{2N}+2\,P_{NG}+P_{2G}>4\,\mu(\tau,c)^2$ (and $P_{2N,e}+2\,P_{NG,e}+P_{2G,e}>4\,\mu_e(\tau)^2$).

\begin{remark}
\textbf{The Joint Distribution of $\left( h_{12},h_{13} \right)$}:
The pair $\left( h_{12},h_{13} \right)$ is a bivariate discrete random
variable with nine possible values such that
$$\left( 2\,h_{12},2\,h_{13} \right)\in \{(0,0),(0,1),(0,2),(1,0),(1,1),(1,2),(2,0),(2,1),(2,2)\}.$$
Then finding the joint distribution of $\left( h_{12},h_{13} \right)$
is equivalent to finding the joint probability mass function of
$\left( h_{12},h_{13} \right)$.
Hence the joint distribution of $\left( h_{12},h_{13} \right)$
can be found by calculating the probabilities such as
$P(\left( h_{12},h_{13} \right)=(0,0))=P(\{X_2,X_3 \} \subset \mI_i \setminus (N(X_1,\tau,c)\cup \G_1(X_1,\tau,c)))$.
$\square$
\end{remark}

\subsection{The Distribution of Relative Density of $\U(y_1,y_2)$-random $\D_{n,2}(\tau,c)$-digraphs}
\label{sec:rel-dens-tau-c}
In the special case of $m=2$ with $\Y_2=\{y_1,y_2\}$ and $\delta_1=y_1<y_2=\delta_2$,
we have only one middle interval and the two end intervals are empty.
In this section, we consider the relative density of central similarity PCD
based on uniform data in $(y_1,y_2)$.
By Theorems \ref{thm:scale-inv-NCSt} and \ref{thm:U(y1,y2-randomDn2(tau,c)},
the asymptotic distribution of any  $\rho_{{}_{[i]}}(\tau,c)$ for the middle intervals for $m>2$
will be identical
to the asymptotic distribution of $\U(y_1,y_2)$-random $\D_{n,2}(\tau,c)$-digraph.

First we consider the simplest case of $\tau=1$ and $c=1/2$.
By Theorem \ref{thm:scale-inv-NCSt},
without loss of generality,
we can assume $(y_1,y_2)$ to be the unit interval $(0,1)$.
Then
$N(x,1,1/2)=B(x,r(x))$ where $r(x)=\min(x,1-x)$ for $x \in (0,1)$.
Hence central similarity PCD based on $N(x,1,1/2)$
is equivalent to the CCCD of \cite{priebe:2001}.
Moreover, we have $\G_1(X_1,2,1/2)=\left(X_1/2,\left(1+X_1\right)/2\right)$.
\begin{theorem}
\label{thm:reldens-tau=1,c=1/2}
As $n \rightarrow \infty$, we have
$\sqrt{n}\,\left[\rho_{n}(1,1/2)-\mu(1,1/2)\right]\stackrel{\mathcal L}{\longrightarrow} \N(0,4\,\nu(1,1/2))$,
where $\mu(1,1/2)=1/2$ and $4\,\nu(1,1/2)=1/12$.
\end{theorem}
\noindent
\textbf{Proof:}
By symmetry, we only consider $X_1 \in (0,1/2)$.
Notice that for $x \in (0,1/2)$,
we have
$N(x,1,1/2)=(0,2\,x)$ and $\G_1(x,1,1/2)=(x/2,(1+x)/2)$.
Hence
$\mu(1,1/2)=P(X_2 \in N(X_1,1,1/2))=2\,P(X_2 \in N(X_1,1,1/2),X_1 \in (0,1/2))$ by symmetry.
Here
\begin{multline*}
P(X_2 \in N(X_1,1,1/2),X_1 \in (0,1/2))=P(X_2 \in (0, 2x_1), X_1 \in (0,1/2))\\
=\int_0^{1/2}\int_0^{2x_1} f_{1,2}(x_1,x_2)dx_2 dx_1
=\int_0^{1/2}\int_0^{2x_1} 1 dx_2 dx_1
=\int_0^{1/2} 2x_1 dx_1= x_1^2|_{0}^{1/2}=1/4.
\end{multline*}
Then $\mu(1,1/2)=2\,(1/4)=1/2$.

For $\Cov(h_{12},h_{13})$, we need to calculate $P_{2N}$, $P_{NG}$, and $P_{2G}$.
The probability
$$
P_{2N}=P(\{X_2,X_3\} \subset N(X_1,1,1/2))
=2\,P(\{X_2,X_3\} \subset N(X_1,1,1/2),X_1 \in (0,1/2))
$$
and
$P(\{X_2,X_3\} \subset N(X_1,1,1/2),X_1 \in (0,1/2))=\int_{0}^{1/2} (2x_1)^2 dx_1=1/6$.
So
$P_{2N}=2\,(1/6)=1/3$.

$P_{NG}=2\,P(X_2 \in N(X_1,1,1/2),X_3 \in \G_1(X_1,1,1/2),X_1 \in (0,1/2))$
and
$$
P(X_2 \in N(X_1,1,1/2),X_3 \in \G_1(X_1,1,1/2),X_1 \in (0,1/2))
=\int_0^{1/2} (2x_1)(1/2)dx_1=1/8.
$$
Then
$P_{NG}=2\,(1/8)=1/4$.

Finally, we have
$P_{2G}=2\,P(\{X_2,X_3\} \subset \G_1(X_1,1,1/2),X_1 \in (0,1/2))$
and
$P(\{X_2,X_3\} \subset \G_1(X_1,1,1/2),X_1 \in (0,1/2))
=\int_0^{1/2} (1/4) dx_1=1/8$.
So
$P_{2G}=2\,(1/8)=1/4$.

Therefore $4\,\E[h_{12}h_{13}]=1/3+2\,(1/4)+1/4=13/12$.
Hence
$4\,\nu(1,1/2)=4\,\Cov[h_{12},h_{13}]=13/12-4(1/2)^2=1/12$.
$\blacksquare$

The sharpest rate of convergence in Theorem \ref{thm:reldens-tau=1,c=1/2} is
$K\,\frac{\mu(2,1/2)}{\sqrt{n\,\nu(2,1/2)^3}}=12 \sqrt{3}\,\frac{K}{\sqrt{n}}$.

Next we consider the more general case of $\tau=1$ and $c \in (0,1)$.
%Without loss of generality,
%assume $0 < c < 1/2$.
For $x \in (0,1)$,
the proximity region has the following form:
\begin{equation}
\label{eqn:NCSt-den-tau=1,c}
N(x,1,c)=
\begin{cases}
(0,x/c) & \text{if $x \in (0,c)$,}\\
((x-c)/(1-c),1) & \text{if $x \in (c,1)$},
\end{cases}
\end{equation}
and the $\G_1$-region is
$\G_1(x,1,c)=(c\,x,(1-c)\,x+c)$.

\begin{theorem}
\label{thm:reldens-tau=1,c}
As $n \rightarrow \infty$,
for $c \in (0,1)$,
we have
$\sqrt{n}\,\left[\rho_{n,2}(1,c)-\mu(1,c)\right]\stackrel{\mathcal L}{\longrightarrow} \N(0,4\,\nu(1,c))$,
where
$\mu(1,c)= 1/2$
and
$4\,\nu(1,c)= c\,(1-c)/3$.
\end{theorem}

Proof is provided in Appendix 1.
See Figure \ref{fig:asy-var-tau=1-c} for $4\,\nu(1,c)$ with $c \in (0,1/2)$.
Notice that $\mu(1,c)$ is constant (i.e., independent of $c$)
and
$\nu(1,c)$ is symmetric around $c=1/2$ with $\nu(1,c)=\nu(1,1-c)$.
Notice also that for $c=1/2$,
we have $\mu(1,c=1/2)=1/2$,
and
$4\,\nu(1,c=1/2)=1/12$,
hence as $c \rightarrow 1/2$,
the distribution of $\rho_{n,2}(1,c)$ converges to the one in Theorem \ref{thm:reldens-tau=1,c=1/2}.
Furthermore,
the sharpest rate of convergence in Theorem \ref{thm:reldens-tau=1,c} is
\begin{equation}
\label{eqn:nu-tau=1,c}
K\,\frac{\mu(1,c)}{\sqrt{n\,\nu(1,c)^3}}=\frac{3\sqrt{3}}{2\,\sqrt{c^3\,(1-c)^3}}\frac{K}{\sqrt{n}}
\end{equation}
and is minimized at $c =1/2$ (which can easily be verified).
% which
%is found by setting the first derivative of this rate
%with respect to $c$ to zero and solving for $c$.
%We also checked the plot of $\mu(1,c)/\sqrt{\nu(1,c)^3}$

\begin{figure}
\centering
\psfrag{c}{\Huge{$c$}}
%\rotatebox{-90}{ \resizebox{2.in}{!}{ \includegraphics{AsyMeanr2c.eps} } }
\rotatebox{-90}{ \resizebox{2.in}{!}{ \includegraphics{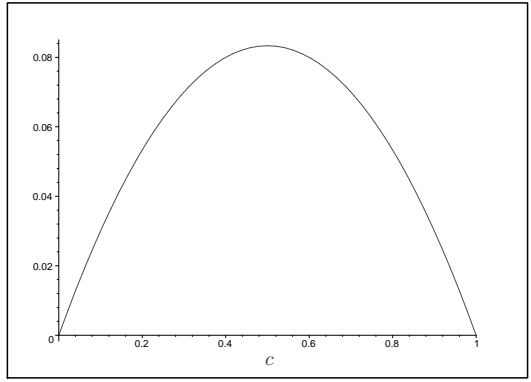} } }
\caption{
\label{fig:asy-var-tau=1-c}
The plot of the asymptotic variance $4\,\nu(1,c)$ as a function of $c$ for $c \in (0,1)$.}
\end{figure}

Next we consider the case of $\tau > 0$ and $c =1/2$.
By symmetry, we only consider $X_1 \in (0,1/2)$.
For $x\in (0,1/2)$,
the proximity region
for $\tau \in (0,1)$ is
\begin{equation}
\label{eqn:NCSt-den-tau01,c=1/2}
N(x,\tau,1/2)=
\begin{cases}
(x\,(1-\tau),x\,(1+\tau)) & \text{if $x \in (0,1/2)$,}\\
(x-(1-x)\,\tau,x+(1-x)\,\tau) & \text{if $x \in (1/2,1)$,}
\end{cases}
\end{equation}
and
for $\tau \ge 1$
\begin{equation}
\label{eqn:NCSt-den-tau-gtr1,c=1/2}
N(x,\tau,1/2)=
\begin{cases}
(0,x\,(1+\tau)) & \text{if $x \in (0,1/(1+\tau))$,}\\
(0,1) & \text{if $x \in (1/(1+\tau),\tau/(1+\tau))$,}\\
(x-(1-x)\,\tau,1) & \text{if $x \in (\tau/(1+\tau),1)$}.
\end{cases}
\end{equation}
And the $\G_1$-region
for $\tau \in (0,1)$ is
\begin{equation}
\label{eqn:GammaCSt-tau-gtr1,c=1/2}
\G_1(x,\tau,1/2)=
\begin{cases}
(x/(1+\tau),x/(1-\tau)) & \text{if $x \in (0,(1-\tau)/2)$,}\\
(x/(1+\tau),(x+\tau)/(1+\tau)) & \text{if $x \in ((1-\tau)/2,(1+\tau)/2)$,}\\
((x-\tau)/(1-\tau),(x+\tau)/(1+\tau)) & \text{if $x \in ((1+\tau)/2,1)$,}
\end{cases}
\end{equation}
and
for $\tau \ge 1$, we have
$\G_1(x,\tau,1/2)=(x/(1+\tau),(x+\tau)/(1+\tau))$.

\begin{theorem}
\label{thm:reldens-tau,c=1/2}
For $\tau \in (0,\infty)$,
we have
$\sqrt{n}\,\left[\rho_{n,2}(\tau,1/2)-\mu(\tau,1/2)\right]\stackrel{\mathcal L}{\longrightarrow} \N(0,4\,\nu(\tau,1/2))$
as $n \rightarrow \infty$,
where
\begin{equation}
\label{eqn:mu-tau,c=1/2}
\mu(\tau,1/2)=
\begin{cases}
\tau/2 & \text{if $0< \tau < 1$,}\\
\tau/(\tau+1) & \text{if $\tau \ge 1$,}
\end{cases}
\end{equation}
and
\begin{equation}
\label{eqn:nu-tau,c=1/2}
4\,\nu(\tau,1/2)=
\begin{cases}
\frac{\tau^2\,(1+2\,\tau-\tau^2-\tau^3)}{3\,(\tau+1)^2} & \text{if $0 < \tau < 1$,}\\
\frac{2\,\tau-1}{3\,(\tau+1)^2}& \text{if $\tau \ge 1$.}
\end{cases}
\end{equation}
\end{theorem}

Proof is provided in Appendix 1.
See Figure \ref{fig:asymean-var-tau-c=1/2} for the plots of $\mu(\tau,1/2)$ and $4\,\nu(\tau,1/2)$.
Notice that $\lim_{\tau \rightarrow \infty} \nu(\tau,1/2)=0$,
so the CLT result fails for $\tau = \infty$.
Furthermore, $\lim_{\tau \rightarrow 0} \nu(\tau,1/2)=0$.
For $\tau=1$,
we have
$\mu(\tau=1,c=1/2)=1/2$,
and
$4\,\nu(\tau=1,c=1/2)=1/12$;
hence as $\tau \rightarrow 1$,
the distribution of $\rho_{n,2}(\tau,1/2)$ converges to the one in Theorem \ref{thm:reldens-tau=1,c=1/2}.
Furthermore,
the sharpest rate of convergence in Theorem \ref{thm:reldens-tau,c=1/2} is
\begin{equation}
\label{eqn:rate-tau,c=1/2}
K\,\frac{\mu(\tau,1/2)}{\sqrt{n\,\nu(\tau,1/2)^3}}=
\frac{K}{\sqrt{n}}\,
\begin{cases}
{\frac {27\,\tau}{2}} \left( {\frac { \left( 6\,\tau+3-3\,\tau^3-3\,\tau^2 \right) \tau^2}
{ \left( \tau+1 \right)^2}}\right)^{-3/2} & \text{if $0 < \tau < 1$,}\\
\frac{3\,\sqrt{3}\,\tau}{\tau+1} \left( {\frac {2\,\tau-1}{ \left( \tau+1 \right)^2}} \right)^{-3/2} & \text{if $\tau \ge 1$.}
\end{cases}
\end{equation}
and is minimized at $\tau \approx .73$ which is found by setting the
first derivative of this rate with respect to $\tau$
to zero and solving for $\tau$ numerically.
We also checked the plot of $\mu(\tau,1/2)/\sqrt{\nu(\tau,1/2)^3}$ (not presented)
and verified that this is where the global minimum is attained.

\begin{figure}
\centering
\psfrag{tau}{\Huge{$\tau$}}
\rotatebox{-90}{ \resizebox{2.in}{!}{ \includegraphics{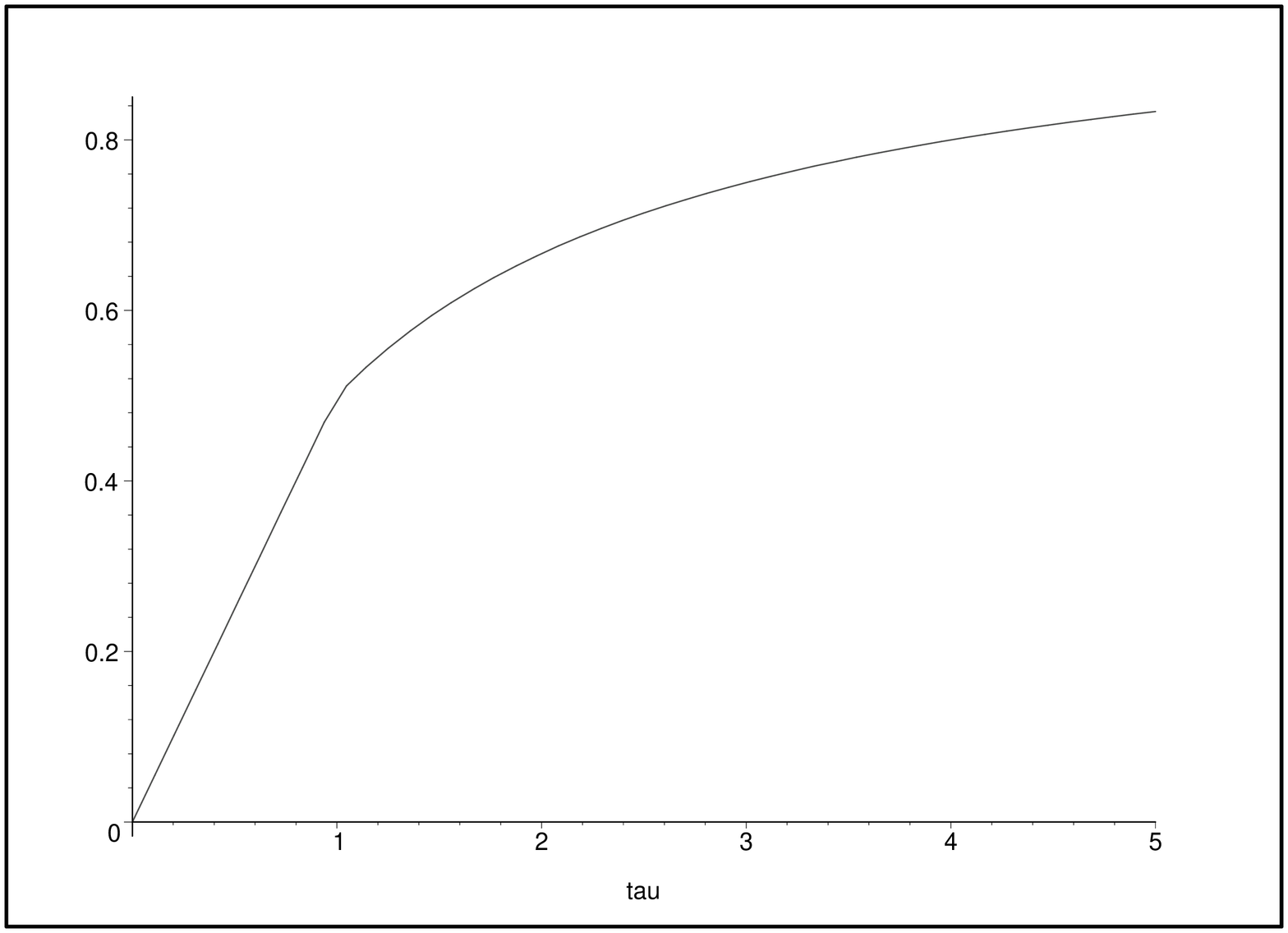} } }
\rotatebox{-90}{ \resizebox{2.in}{!}{ \includegraphics{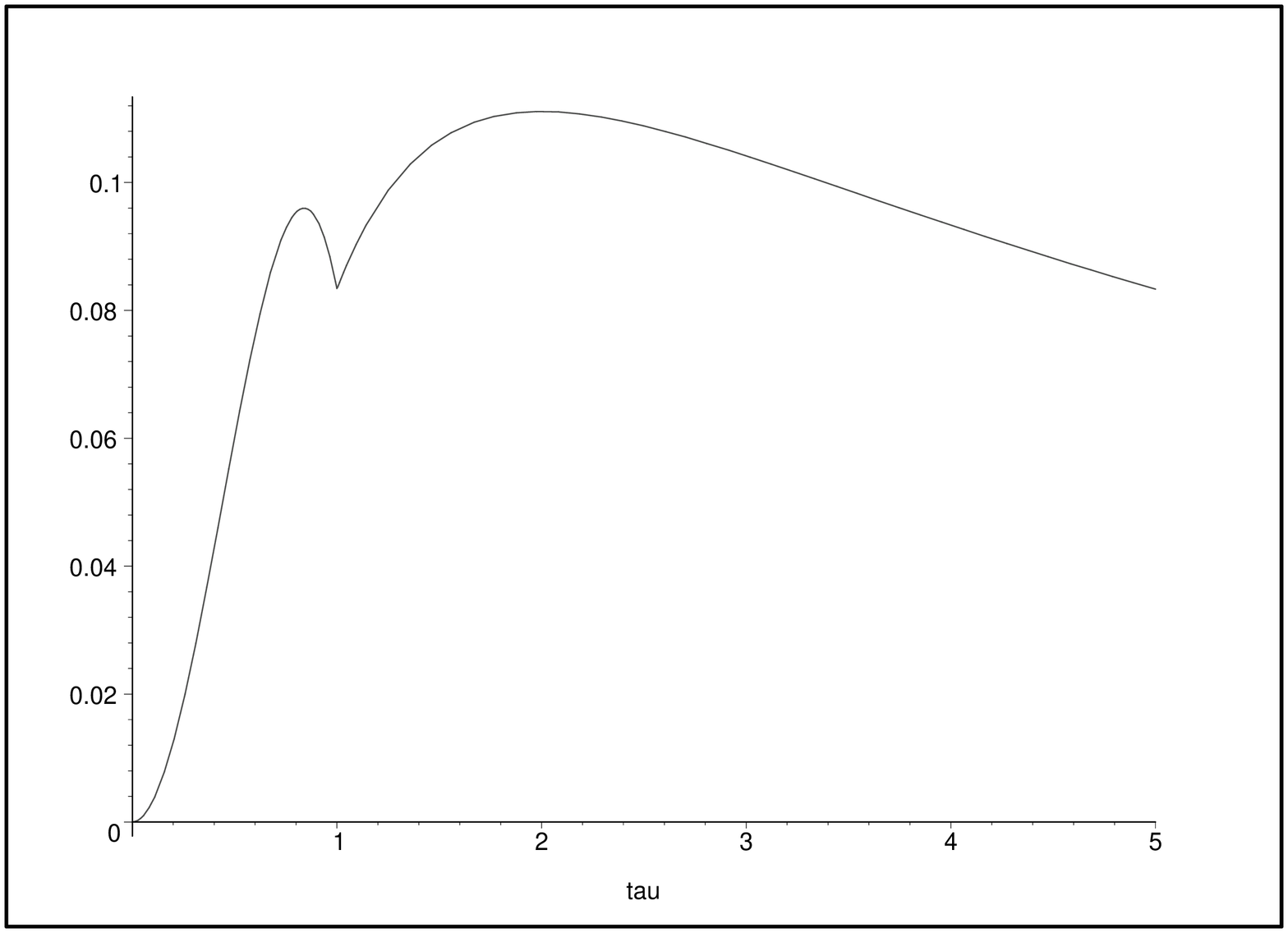} } }
\caption{
\label{fig:asymean-var-tau-c=1/2}
The plots of the asymptotic mean $\mu(\tau,1/2)$ (left) and the variance $4\,\nu(\tau,1/2)$ (right)
as a function of $\tau$ for $\tau \in (0,5]$.}
\end{figure}

Finally,
we consider the most general case of
$\tau > 0$ and $c \in (0,1/2)$.
For $\tau \in (0,1)$,
the proximity region is
\begin{equation}
\label{eqn:NCS-defn-tau01}
N(x,\tau,c)=
\begin{cases}
\left(x\,(1-\tau),x\,\left(1+\frac{(1-c)\,\tau}{c}\right)\right) & \text{if $x \in (0,c)$,} \vspace{.1cm}\\
\left(x-\frac{c\,\tau\,(1-x)}{1-c},x+(1-x)\,\tau\right) & \text{if $x \in (c,1)$,}
\end{cases}
\end{equation}
and the $\G_1$-region is
\begin{equation}
\label{eqn:G1-defn-tau01}
\G_1(x,\tau,c)=
\begin{cases}
\left(\frac{c\,x}{c+(1-c)\,\tau},\frac{x}{1-\tau}\right) & \text{if $x \in (0,c\,(1-\tau))$,}\vspace{.1cm}\\
\left(\frac{c\,x}{c+(1-c)\,\tau},\frac{x\,(1-c)+c\,\tau}{1-c+c\,\tau}\right) & \text{if $x \in (c\,(1-\tau),c\,(1-\tau)+\tau)$,}\vspace{.1cm}\\
\left(\frac{x-\tau}{1-\tau},\frac{x\,(1-c)+c\,\tau}{1-c+c\,\tau}\right) & \text{if $x \in (c\,(1-\tau)+\tau,1)$.}
\end{cases}
\end{equation}

For $\tau \ge 1$,
the proximity region is
\begin{equation}
\label{eqn:NCS-defn-taugtr1}
N(x,\tau,c)=
\begin{cases}
\left(0,x\,\left(1+\frac{(1-c)\,\tau}{c}\right)\right) & \text{if $x \in \left(0,\frac{c}{c+(1-c)\,\tau}\right)$,}\vspace{.1cm}\\
\left(0,1\right) & \text{if $x \in \left(\frac{c}{c+(1-c)\,\tau},\frac{c\,\tau}{1-c+c\,\tau}\right)$,}\vspace{.1cm}\\
\left(x-\frac{c\,\tau\,(1-x)}{1-c},1\right) & \text{if $x \in \left(\frac{c\,\tau}{1-c+c\,\tau},1\right)$,}
\end{cases}
\end{equation}
and the $\G_1$-region is
\begin{equation}
\label{eqn:G1-defn-taugtr1}
\G_1(x,\tau,c)=
\left(\frac{c\,x}{c+(1-c)\,\tau},\frac{x\,(1-c)+c\,\tau}{1-c+c\,\tau}\right).
\end{equation}

\begin{theorem}
\label{thm:reldens-tau,c}
For $\tau \in (0,\infty)$,
we have
$\sqrt{n}\,\left[\rho_{n,2}(\tau,c)-\mu(\tau,c)\right]\stackrel{\mathcal L}{\longrightarrow} \N(0,4\,\nu(\tau,c))$,
as $n \rightarrow \infty$,
where
$\mu(\tau,c)=\mu_1(\tau,c)\,\I(0 < c \le 1/2)$ + $\mu_2(\tau,c)\,\I(1/2 \le c < 1)$
and
$\nu(\tau,c)=\nu_1(\tau,c)\,\I(0 < c \le 1/2)+\nu_2(\tau,c)\,\I(1/2 \le c < 1)$.
For $0 < c \le 1/2$,
\begin{equation}
\label{eqn:mu-tau,c}
\mu_1(\tau,c)=
\begin{cases}
\frac{\tau}{2} & \text{if $0< \tau < 1$,}\\
\frac{\tau\,(1+2\,c\,(\tau-1)(1-c))}{2\,(c\,\tau-c+1)(\tau+c-c\,\tau)} & \text{if $\tau \ge 1$,}
\end{cases}
\end{equation}
and
\begin{equation}
\label{eqn:nu-tau,c}
4\,\nu_1(\tau,c)=
\begin{cases}
\kappa_1(\tau,c) & \text{if $0< \tau < 1$,}\\
\kappa_2(\tau,c) & \text{if $\tau \ge 1$,}
\end{cases}
\end{equation}
where
$$
\kappa_1(\tau,c)= \frac {\tau^2 \left( c^2\,\tau^3-3\,c^2\,\tau^2
-c\,\tau^3+2\,c^2\,\tau+3\,c\,\tau^2-c^2-2\,c\,\tau-\tau^2
+c+\tau \right) }{ 3\,\left( c\,\tau-c+1 \right)  \left( c+\tau-c\,\tau\right) },
$$
and
{\small
\begin{multline*}
\kappa_2(\tau,c)=\Bigl[c (1-c)  \bigl( 2\,c^4\,\tau^5-7\,c^4\,\tau^4-4\,c^3\,\tau^5+
8\,c^4\,\tau^3+14\,c^3\,\tau^4+3\,c^2\,\tau^5-2\,c^4\,\tau^2-16\,c^3\,\tau^3-7\,c^2\,\tau^4-c\,\tau^5-\\
2\,c^4\,\tau+4\,c^3\,\tau^2+12\,c^2\,\tau^3+c^4+4\,c^3\,\tau-6\,c^2\,\tau^2-
4\,c\,\tau^3-2\,c^3-3\,c^2\,\tau+4\,c\,\tau^2+c^2+c\,\tau-\tau^2 \bigr)\Bigr]\Big/
\Bigl[3\, \left( c\,\tau-c+1 \right)^3 \left( c\,\tau-c-\tau \right)^3\Bigr].
\end{multline*}
}
And for $1/2 \le c < 1$,
we have
$\mu_2(\tau,c)=\mu_1(\tau,1-c)$
and
$\nu_2(\tau,c)=\nu_1(\tau,1-c)$.
\end{theorem}

Proof is provided in Appendix 1.
See Figure \ref{fig:asymean-var-tau-c} for the plots of
$\mu(\tau,c)$ and $4\,\nu(\tau,c)$.
Notice that $\lim_{\tau \rightarrow \infty} \nu(\tau,c)=0$,
so the CLT result fails for $\tau=\infty$.
Furthermore,
$\lim_{\tau \rightarrow 0}\nu(\tau,c)=0$.
For $\tau=1$ and $c=1/2$,
we have
$\mu(\tau=1,c=1/2)=1/2$,
and
$4\,\nu(\tau=1,c=1/2)=1/12$,
hence as $\tau \rightarrow 1$ and $c \rightarrow 1/2$,
the distribution of $\rho_{n,2}(\tau,c)$ converges to the one in Theorem \ref{thm:reldens-tau=1,c=1/2}.
The sharpest rate of convergence in Theorem \ref{thm:reldens-tau,c} is
$K\,\frac{\mu(\tau,c)}{\sqrt{n\,\nu(\tau,c)^3}}$ (the explicit form not presented)
and is minimized at $\tau \approx 1.55$ and $c \approx 0.5$
which is found by setting the first order partial derivatives
of this rate with respect to $\tau$ and $c$ to zero and
solving for $\tau$ and $c$ numerically.
We also checked the surface plot of this rate (not presented)
and verified that this is where the global minimum is attained.

\begin{figure}
\centering
\psfrag{t}{\Huge{$\tau$}}
\psfrag{c}{\Huge{$c$}}
\rotatebox{0}{ \resizebox{2.9in}{!}{ \includegraphics{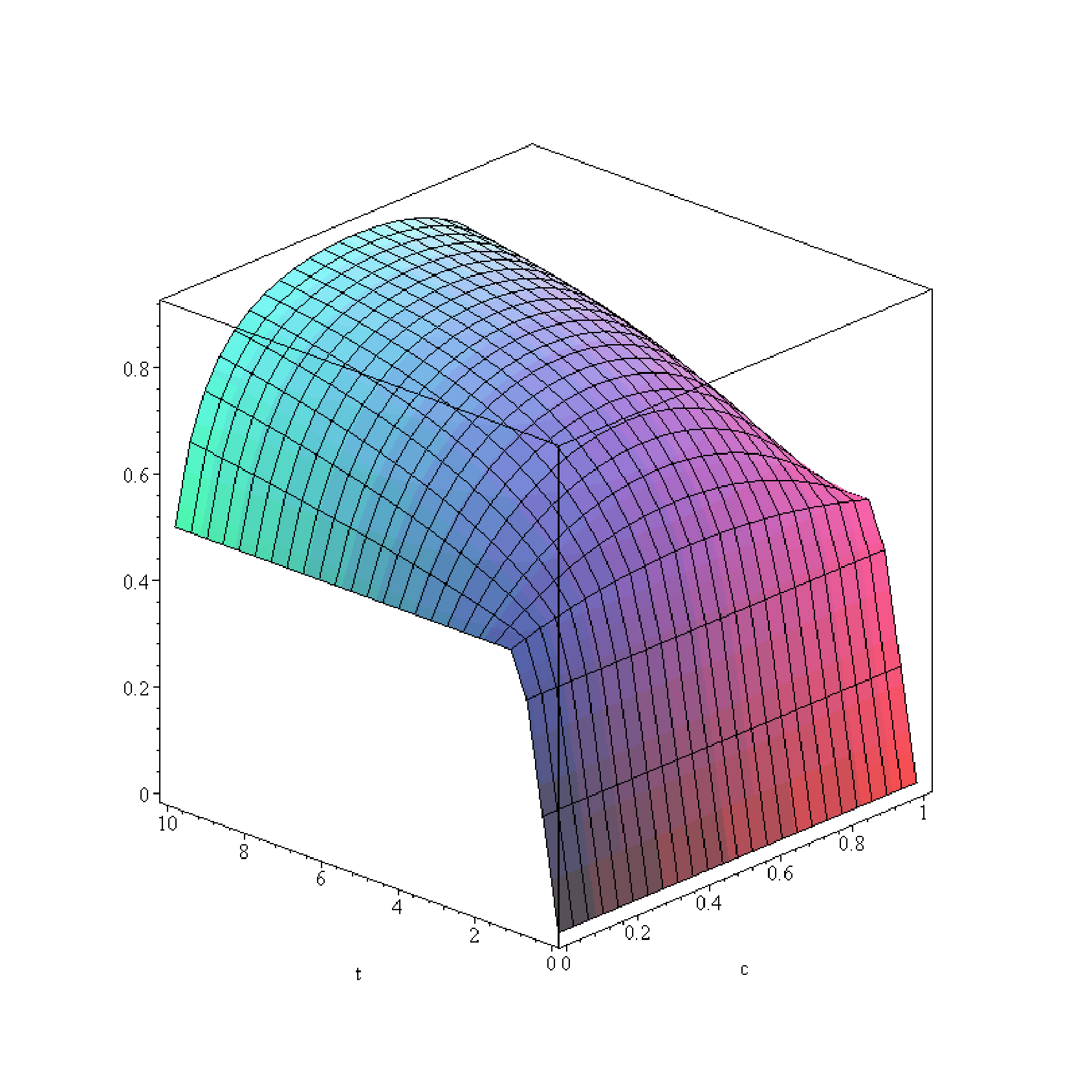} } }
\rotatebox{0}{ \resizebox{2.9in}{!}{ \includegraphics{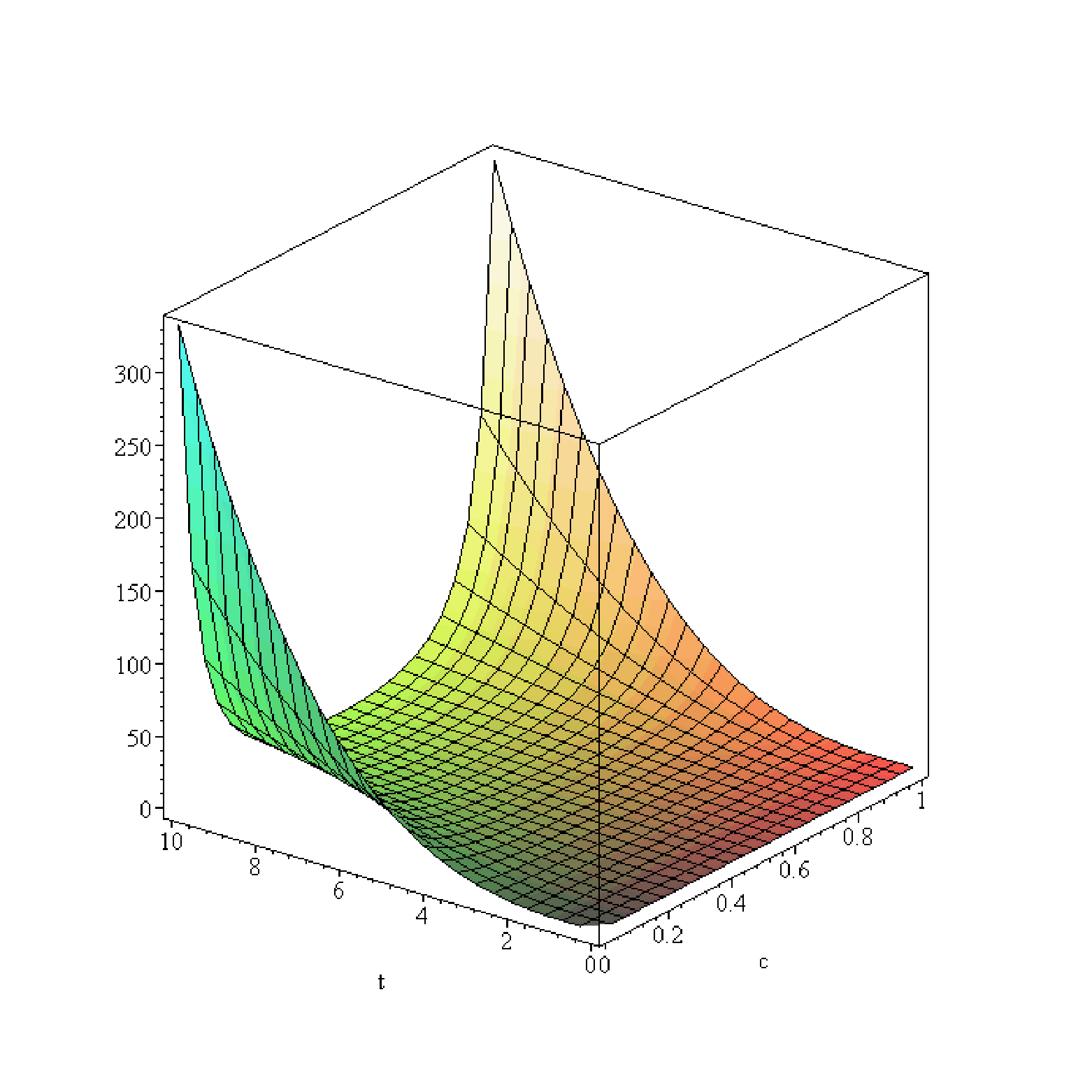} } }
\caption{
\label{fig:asymean-var-tau-c}
The surface plots of the asymptotic mean $\mu(t,c)$ (left) and the variance $4\,\nu(t,c)$ (right)
as a function of $t$ and $c$ for $t \in (0,10]$ and $c \in (0,1)$, respectively.}
\end{figure}

\subsection{The Case of End Intervals:
Relative Density for $\U\left(\delta_1,y_{(1)}\right)$ or $\U\left(y_{(m)},\delta_2\right)$ Data}
\label{sec:rel-dens-end-intervals}
Recall that with $m \ge 1$ for the end intervals,
$\mI_1=\left(\delta_1,y_{(1)}\right)$ and $\mI_{m+1}=\left(y_{(m)},\delta_2\right)$,
the proximity and $\G_1$-regions were only dependent on $x$
and $\tau$ (but not on $c$).
%Hence the new notation $\D_{n,2}(\tau)$-digraph in the title.
%For $\U\left(\delta_1,y_{(1)}\right)$ and $\U\left(y_{(m)},\delta_2\right)$ data,
%by symmetry, relative density has the same distribution.
%So we only consider $\left(y_{(m)},\delta_2\right)$.
Due to scale invariance from Theorem \ref{thm:scale-inv-NCSt},
we can assume that each of the end intervals is $(0,1)$.
Let %$N_{e}(x,\tau)$ be the central similarity proximity region,
$\G_{1,e}(x,\tau)$ be the $\G_1$-region corresponding to $N_{e}(x,\tau)$ in the end interval case.
%Moreover,
%let $\mu_e(\tau)$ and $4\,\nu_e(\tau)$ be the asymptotic mean and variance
%for the relative density of the $\D_{n,2}(\tau)$-digraph in the end interval case.

First we consider $\tau=1$ and uniform data in the end intervals.
Then for $x$ in the right end interval,
$N_{e}(x,1)=(0,\min(1,2x))$ for $x \in (0,1)$
and the $\G_1$-region is
$\G_{1,e}(x,1)=(x/2,1)$.
\begin{theorem}
\label{thm:reldens-end-int-tau=1}
Let $D_{[i]}(1,c)$ be the subdigraph of the central similarity PCD based on uniform data in $(\delta_1,\delta_2)$
where $\delta_1<\delta_2$ and $\Y_m$ be a set of $m$ distinct $\Y$ points in $(\delta_1,\delta_2)$.
Then for $i \in \{1,m+1\}$ (i.e., in the end intervals),
as $n_i \rightarrow \infty$, we have
$\sqrt{n_i}\,\left[\rho_{{}_{[i]}}(1,c)-\mu_e(1)\right]\stackrel{\mathcal L}{\longrightarrow} \N(0,4\,\nu_e(1))$,
where $\mu_e(1)=3/4$ and $4\,\nu_e(1)=1/24$.
\end{theorem}

The Proof is provided in Appendix 1.
The sharpest rate of convergence in Theorem \ref{thm:reldens-end-int-tau=1} is
$K\,\frac{\mu_e(1)}{\sqrt{n_i\,\nu_e(1)^3}}=36 \sqrt{6}\,\frac{K}{\sqrt{n_i}}$ for $i \in \{1,m+1\}$.

Next we consider the more general case of $\tau>0$ for the end intervals.
By Theorem \ref{thm:scale-inv-NCSt},
we can assume each end interval to be $(0,1)$.
For $\tau \in (0,1)$
and $x$ in the right end interval,
the proximity region is
\begin{equation}
\label{eqn:NCS-defn-tau01-end}
N_e(x,\tau)=
\begin{cases}
\left(x\,(1-\tau),x\,(1+\tau)\right) & \text{if $x \in (0,1/(1+\tau))$,}\\
\left(x\,(1-\tau),1\right) & \text{if $x \in (1/(1+\tau),1)$,}\\
\end{cases}
\end{equation}
and the $\G_1$-region is
\begin{equation}
\label{eqn:G1-defn-tau01-end}
\G_{1,e}(x,\tau)=
\begin{cases}
\left(\frac{x}{1+\tau},\frac{x}{1-\tau}\right) & \text{if $x \in (0,1-\tau)$,}\vspace{.1cm}\\
\left(\frac{x}{1+\tau},1\right) & \text{if $x \in (1-\tau,1)$.}
\end{cases}
\end{equation}

For $\tau \ge 1$
and $x$ in the right end interval,
the proximity region is
\begin{equation}
\label{eqn:NCS-defn-taugtr1-end}
N_e(x,\tau)=
\begin{cases}
\left(0,x\,(1+\tau)\right) & \text{if $x \in (0,1/(1+\tau))$,}\\
\left(0,1\right) & \text{if $x \in (1/(1+\tau),1)$,}\\
\end{cases}
\end{equation}
and the $\G_1$-region is
$\G_{1,e}(x,\tau)=\left(x/(1+\tau),1\right).$

\begin{theorem}
\label{thm:reldens-end-int-tau}
Let $D_{[i]}(\tau,c)$ be the subdigraph of the central similarity PCD based on uniform data in $(\delta_1,\delta_2)$
where $\delta_1<\delta_2$ and $\Y_m$ be a set of $m$ distinct $\Y$ points in $(\delta_1,\delta_2)$.
Then for $i \in \{1,m+1\}$ (i.e., in the end intervals),
and
$\tau \in (0,\infty)$,
we have $\sqrt{n_i}\,\left[\rho_{{}_{[i]}}(\tau,c)-\mu_e(\tau)\right]\stackrel{\mathcal L}{\longrightarrow} \N(0,4\,\nu_e(\tau))$,
as $n_i \rightarrow \infty$,
where
\begin{equation}
\label{eqn:mu-tau,c-end}
\mu_e(\tau)=
\begin{cases}
\frac{\tau\,(\tau+2)}{2\,(\tau+1)} & \text{if $0< \tau < 1$,}\\
\frac{1+2\,\tau}{2\,(\tau+1)} & \text{if $\tau \ge 1$,}
\end{cases}
\end{equation}
and
\begin{equation}
\label{eqn:nu-tau,c-end}
4\,\nu_e(\tau)=
\begin{cases}
{\frac {\tau^2 \left( 4\,\tau+4-2\,\tau^4-4\,\tau^3-\tau^2 \right) }
{ 3\,\left( \tau+1 \right)^3}} & \text{if $0< \tau < 1$,}\\
\frac{\tau^2}{3\,(\tau+1)^3} & \text{if $\tau \ge 1$.}
\end{cases}
\end{equation}
\end{theorem}

See Appendix 1 for the proof
and Figure \ref{fig:asymean-var-tau-end} for the plots of $\mu_e(\tau)$ and $4\,\nu_e(\tau)$.
Notice that $\lim_{\tau \rightarrow \infty} \nu_e(\tau)=0$,
so the CLT result fails for $\tau =\infty$.
Furthermore,
$\lim_{\tau \rightarrow 0}\nu_e(\tau)=0$.
For $\tau=1$, we have $\mu_e(\tau=1)=3/4$,
and
$4\,\nu_e(\tau=1)=1/24$,
hence as $\tau \rightarrow 1$,
the distribution of $\rho_{{}_{[i]}}(\tau,c)$ converges to the one in Theorem \ref{thm:reldens-end-int-tau=1} for $i \in \{1,m+1\}$..
The sharpest rate of convergence in Theorem \ref{thm:reldens-end-int-tau} is
$K\,\frac{\mu_e(\tau)}{\sqrt{n_i\,\nu_e(\tau)^3}}$ (explicit form not presented) for $i \in \{1,m+1\}$
and is minimized at $\tau \approx 0.58$ which is found numerically as before.
We also checked the plot of $\mu_e(\tau)/\sqrt{\nu_e(\tau)^3}$ (not presented)
and verified that this is where the global minimum is attained.

\begin{figure}
\centering
\psfrag{tau}{\Huge{$\tau$}}
\rotatebox{-90}{ \resizebox{2.in}{!}{ \includegraphics{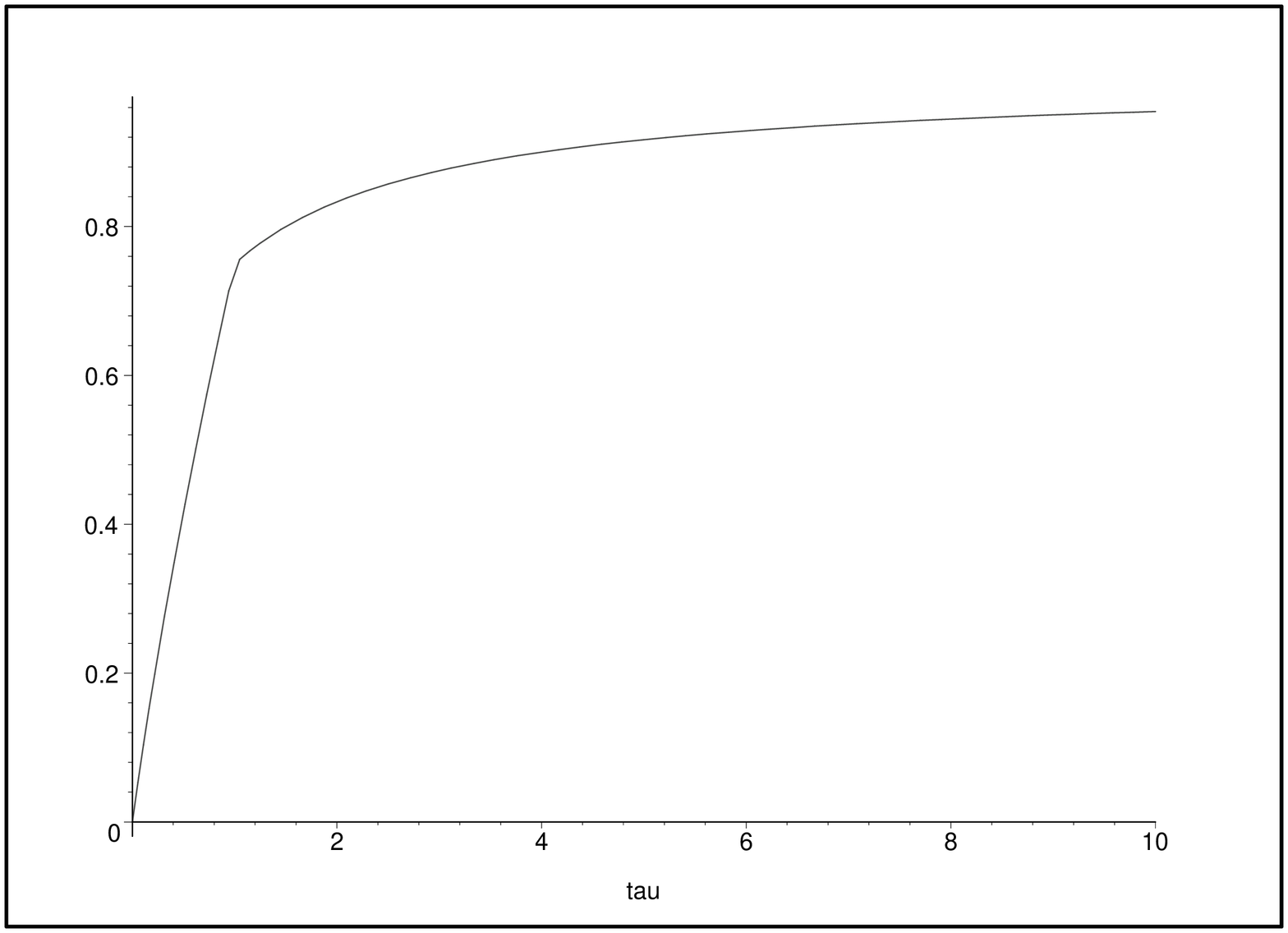} } }
\rotatebox{-90}{ \resizebox{2.in}{!}{ \includegraphics{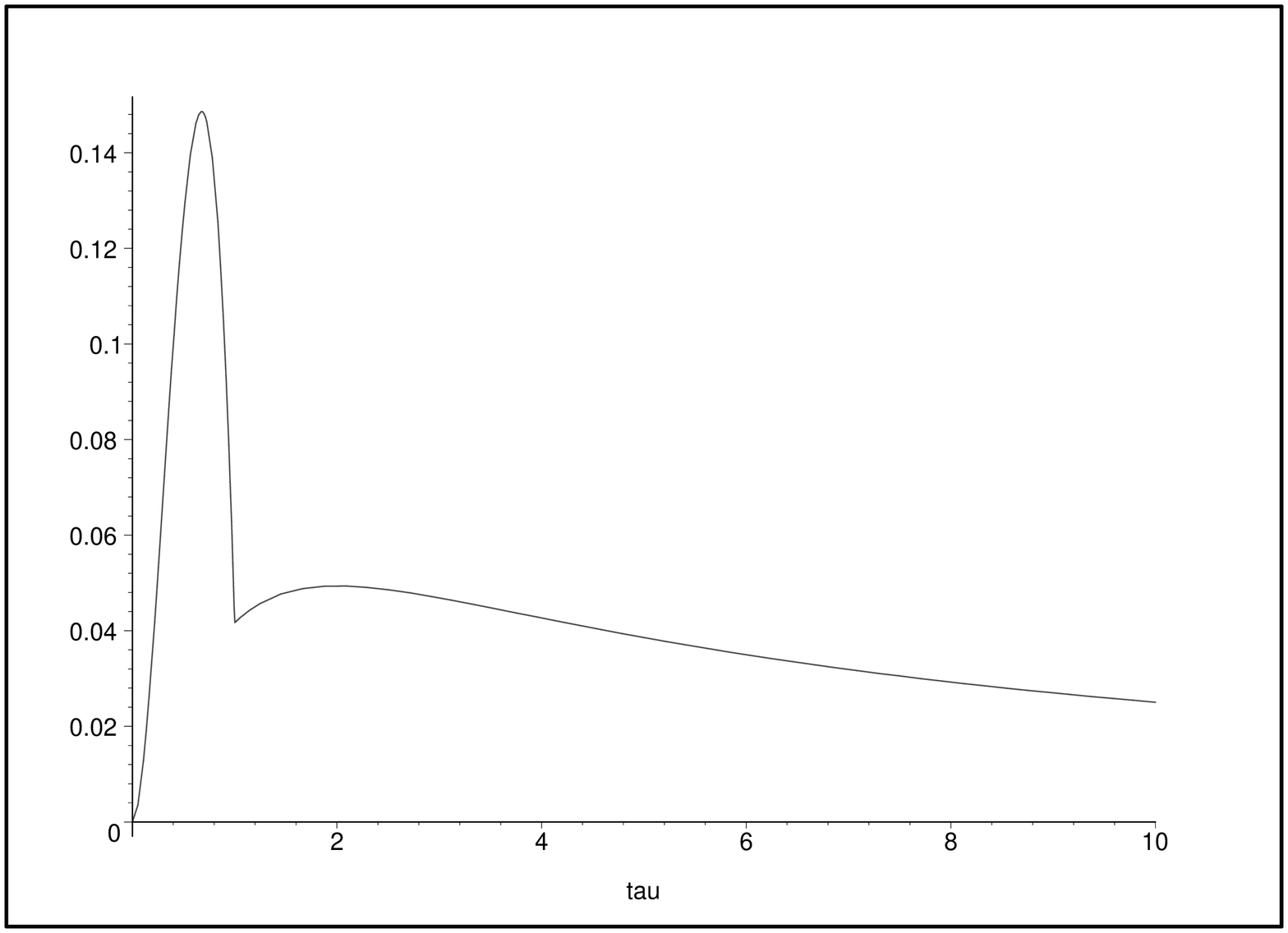} } }
\caption{
\label{fig:asymean-var-tau-end}
The plots of the asymptotic mean $\mu_e(\tau)$ (left) and the variance $4\,\nu_e(\tau)$ (right)
for the end intervals as a function of $\tau$ for $\tau \in (0,10]$.}
\end{figure}

\section{The Distribution of the Relative Density of $\U(\delta_1,\delta_2)$-random $\D_{n,m}(\tau,c)$-digraphs}
\label{sec:dist-multiple-intervals}
In this section,
we consider the more challenging case of $m>2$.

\subsection{First Version of Relative Density in the Case of $m \ge 2$}
\label{sec:version-I-mult-int}
Recall that the relative density $\displaystyle \rho_{n,m}(\tau,c)$
is defined as in Equation \eqref{eqn:rho-nm-tau,c}.
%So $\displaystyle \rho_{n,m}(\tau,c)=\frac{\sum_{i=1}^{m+1} |\A_{[i]}|}{n_{_T}}=
%\sum_{i=1}^{m+1}\frac{n_i\,(n_i-1)}{n_{_T}}\,\rho_{{}_{[i]}}(\tau,c)$.
%Since $\frac{n_i\,(n_i-1)}{n_{_T}} \ge 0$ for each $i$
%and $\displaystyle \sum_{i=1}^{m+1}\frac{n_i\,(n_i-1)}{n_{_T}}=1$,
%it follows that
%$\rho_{n,m}(\tau,c)$ is a mixture of the $\rho_{{}_{[i]}}(\tau,c)$.
Letting $w_i = \left(y_{(i+1)}-y_{(i)}\right) / (\delta_2-\delta_1)$, for $i=0,1,2,\ldots,m$,
we obtain the following as a result of Theorem \ref{thm:reldens-tau,c}.

\begin{theorem}
\label{thm:MI-asy-norm-I}
Let $\X_n$ be a random sample from $\U(\delta_1,\delta_2)$
with $-\infty<\delta_1<\delta_2<\infty$ and $\Y_m$ be a set of $m$ distinct points in $(\delta_1,\delta_2)$.
For $\tau \in (0,\infty)$,
the asymptotic distribution of $\rho_{n,m}(\tau,c)$ conditional on $\Y_m$
is given by
\begin{equation}
\sqrt{n}\left(\rho_{n,m}(\tau,c)-\breve \mu(m,\tau,c)\right)
\stackrel{\mathcal L}{\longrightarrow}\\
\mathcal{N}
 \left(
   0,
   4\,\breve \nu(m,\tau,c)
 \right),
\end{equation}
as $n \rightarrow \infty$,
provided that $\breve \nu(m,\tau,c)>0$,
where
$\breve \mu(m,\tau,c)=\widetilde \mu(m,\tau,c)\Big/\left(\sum_{i=1}^{m+1} w_i^2 \right)$
with
$\widetilde \mu(m,\tau,c)=\mu(\tau,c)\sum_{i=2}^{m} w_i^2+\mu_e(\tau)\sum_{i \in \{1,m+1\}} w_i^2$
and $\mu(\tau,c)$ and $\mu_e(\tau)$ are as in Theorems \ref{thm:reldens-tau,c} and \ref{thm:reldens-end-int-tau},
respectively.
Furthermore,
$4 \breve \nu(m,\tau,c)= 4 \widetilde \nu(m,\tau,c)\Big/\left(\sum_{i=1}^{m+1} w_i^2 \right)^2$
with
$4\widetilde \nu(m,\tau,c)=
[P_{2N}+2\,P_{NG}+P_{2G}]\sum_{i=2}^{m} w_i^3+[P_{2N,e}+2\,P_{NG,e}+P_{2G,e}]\sum_{i \in \{1,m+1\}} w_i^3- (\widetilde \mu(m,\tau,c))^2$.
%with $P_{2N},2\,P_{NG},P_{2G}$ being as in the Proof of Theorem \ref{thm:reldens-tau,c}
%and
%$P_{2N,e},2\,P_{NG,e},P_{2G,e}$ being as in the Proof of Theorem \ref{thm:reldens-end-int-tau}.
\end{theorem}

Proof is provided in Appendix 2.
Notice that if $y_{(1)}=\delta_1$ and $y_{(m)}=\delta_2$,
there are only $m-1$ middle intervals formed by $y_{(i)}$.
That is, the end intervals $\mI_{1}=\mI_{m+1}=\emptyset$.
Hence in Theorem \ref{thm:MI-asy-norm-I},
$\breve \mu(m,\tau,c)=\mu(\tau,c)$
since
$\widetilde \mu(m,\tau,c)=\mu(\tau,c)\sum_{i=2}^{m} w_i^2$.
Furthermore,
$4 \breve \nu(m,\tau,c)=
[P_{2N}+2\,P_{NG}+P_{2G}]\sum_{i=2}^{m} w_i^3- (\mu(\tau,c)\sum_{i=2}^{m} w_i^2)^2=
4\,\nu(m,\tau,c)+\mu^2(\tau,c)\left(\sum_{i=2}^{m} w_i^3-(\sum_{i=2}^{m} w_i^2)^2\right)$.

\subsection{Second Version of Relative Density in the Case of $m \ge 2$}
\label{sec:version-II-mult-int}
For $m \ge 2$,
if we consider the entire data set $\X_n$,
then we have $n$ vertices.
So we can also consider the relative density as
$\widetilde{\rho}_{n,m}(\tau,c)=\left|\A\right|/(n\,(n-1))$.

\begin{theorem}
\label{thm:MI-asy-norm-II}
Let $\X_n$ be a random sample from $\U(\delta_1,\delta_2)$
with $-\infty<\delta_1<\delta_2<\infty$ and $\Y_m$ be a set of $m$ distinct points in $(\delta_1,\delta_2)$.
For $\tau \in (0,\infty)$,
the asymptotic distribution for $\widetilde{\rho}_{n,m}(\tau,c)$ conditional on $\Y_m$
is given by
\begin{equation}
\sqrt{n}\left(\widetilde{\rho}_{n,m}(\tau,c)-\widetilde \mu(m,\tau,c)\right)
\stackrel{\mathcal L}{\longrightarrow}\\
\mathcal{N}
 \left(
   0,
   4\,\widetilde \nu(m,\tau,c)
 \right),
\end{equation}
as $n \rightarrow \infty$, provided that $\widetilde \nu(m,\tau,c)>0$,
where $\widetilde \mu(m,\tau,c)$
and $\widetilde \nu(m,\tau,c)$
are as in Theorem \ref{thm:MI-asy-norm-I}.
\end{theorem}

Proof is provided in Appendix 2.
Notice that the relative arc densities,
$\rho_{n,m}(\tau,c)$ and $\widetilde \rho_{n,m}(\tau,c)$
do not have the same distribution for neither finite nor infinite $n$.
But we have $\rho_{n,m}(\tau,c)=\frac{n(n-1)}{n_{_T}} \widetilde \rho_{n,m}(\tau,c)$
and
since for large $n_i$ and $n$,
$\sum_{i=1}^{m+1} \frac{n_i(n_i-1)}{n(n-1)} \approx \sum_{i=1}^{m+1} w_i^2<1$,
it follows that
$\widetilde \mu(m,\tau,c) < \breve \mu(m,\tau,c)$
and
$\widetilde \nu(m,\tau,c) < \breve \nu(m,\tau,c)$ for large $n_i$ and $n$.
Furthermore,
the asymptotic normality holds for $\rho_{n,m}(\tau,c)$
iff
it holds for $\widetilde \rho_{n,m}(\tau,c)$.

\section{Extension of Central Similarity Proximity Regions to Higher Dimensions}
\label{sec:NCSt-higher-D}
Note that in $\R$ the central similarity PCDs are based on the intervals
whose end points are from class $\Y$.
This interval partitioning can
be viewed as the \emph{Delaunay tessellation} of $\R$ based on $\Y_m$. So
in higher dimensions, we use the Delaunay tessellation based
on $\Y_m$ to partition the space.

Let $\Y_m=\left \{\y_1,\y_2,\ldots,\y_m \right\}$ be $m$ points in
general position in $\R^d$ and $T_i$ be the $i^{th}$ Delaunay cell
for $i=1,2,\ldots,J_m$, where $J_m$ is the number of Delaunay cells.
Let $\X_n$ be a set of iid random variables from distribution $F$ in
$\R^d$ with support $\mS(F) \subseteq \C_H(\Y_m)$
where $\C_H(\Y_m)$ stands for the convex hull of $\Y_m$.

\subsection{Extension of Central Similarity Proximity Regions to $\mathbb R^2$}
\label{sec:extension-PEPCD-R2}
For illustrative purposes, we focus on $\R^2$ where
a Delaunay tessellation is a \emph{triangulation}, provided that no more
than three points in $\Y_m$ are cocircular (i.e., lie on the same circle).
Furthermore, for simplicity, we only consider the one Delaunay triangle case.
Let $\Y_3=\{\y_1,\y_2,\y_3\}$ be three non-collinear points
in $\R^2$ and $\TY=T(\y_1,\y_2,\y_3)$ be the triangle
with vertices $\Y_3$.
Let $\X_n$ be a set of iid random variables from $F$ with
support $\mS(F) \subseteq \TY$.
%If $F=\U(\TY)$, a composition of translation,
%rotation, reflections, and scaling
%will take any given triangle $\TY$
%to the basic triangle $T_b=T((0,0),(1,0),(c_1,c_2))$
%with $0 < c_1 \le 1/2$, $c_2>0$,
%and $(1-c_1)^2+c_2^2 \le 1$, preserving uniformity.
%That is, if $X \sim \U(\TY)$ is transformed in the same manner to,
%say $X'$, then we have $X' \sim \U(T_b)$.
%In fact this will hold for any distribution $F$
%up to scale.

\begin{figure} [ht]
\centering
    \scalebox{.35}{\input{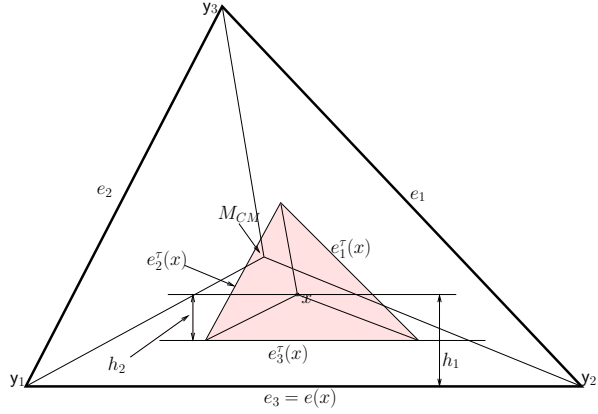}}
\caption{
\label{fig:prox-map-def}
Construction of central similarity proximity region, $N(x,\tau=1/2,M_C)$ (shaded region)
for an $x \in R_E(e_3)$
where $h_2=d(x,e_3^\tau(x))=\frac{1}{2}\,d(x,e(x))$ and $h_1=d(x,e(x))$..}
\end{figure}

\label{sec:tau-factor-PCD}
For the expansion parameter $\tau \in (0,\infty]$,
define $N(x,\tau,M_C)$ to be the {\em central similarity proximity map} with expansion parameter $\tau$ as follows;
see also Figure \ref{fig:prox-map-def}.
Let $e_j$ be the edge opposite vertex $\y_j$ for $j=1,2,3$,
and let ``edge regions'' $R_E(e_1)$, $R_E(e_2)$, $R_E(e_3)$
partition $\TY$ using line segments from the
center of mass of $\TY$ to the vertices.
For $x \in (\TY)^o$, let $e(x)$ be the
edge in whose region $x$ falls; $x \in R_E(e(x))$.
If $x$ falls on the boundary of two edge regions we assign $e(x)$ arbitrarily.
For $\tau > 0$, the central similarity proximity region
$N(x,\tau,M_C)$ is defined to be the triangle $T_{CS}(x,\tau) \cap \TY$ with the following properties:
\begin{itemize}
\item[(i)]
For $\tau \in (0,1]$,
the triangle
$T_{CS}(x,\tau)$ has an edge $e_\tau(x)$ parallel to $e(x)$ such that
$d(x,e_\tau(x))=\tau\, d(x,e(x))$
and
$d(e_\tau(x),e(x)) \le d(x,e(x))$
and
for $\tau >1$,
$d(e_\tau(x),e(x)) < d(x,e_\tau(x))$
where $d(x,e(x))$ is the Euclidean distance from $x$ to $e(x)$,
\item[(ii)] the triangle $T_{CS}(x,\tau)$ has the same orientation as and is similar to $\TY$,
\item[(iii)] the point $x$ is at the center of mass of $T_{CS}(x,\tau)$.
\end{itemize}
Note that (i) implies the expansion parameter $\tau$,
(ii) implies ``similarity", and
(iii) implies ``central" in the name, (parameterized) {\em central similarity proximity map}.
Notice that $\tau>0$ implies that $x \in N(x,\tau,M_C)$ and,
by construction, we have
$N(x,\tau,M_C)\subseteq \TY$ for all $x \in \TY$.
For $x \in \partial(\TY)$ and $\tau \in (0,\infty]$, we define $N(x,\tau,M_C)=\{x\}$.
%and for $\tau=0$ and $x \in \TY$ we define $N(x,\tau,M_C)=\{x\}$.
For all $ x\in \TY^o$ the edges
$e_\tau(x)$ and $e(x)$ are coincident iff $\tau=1$.
Note also that
$\lim_{\tau \rightarrow \infty} N(x,\tau,M_C) = \TY$
for all $x \in (\TY)^o$,
so we define $N(x,\infty,M_C) = \TY$ for all such $x$.

\subsection{Extension of Central Similarity Proximity Regions to $\mathbb R^d$ with $d>2$}
\label{sec:extension-PEPCD-Rd}
The extension to $\R^d$ for $d > 2$ with
$M=M_C$ is provided in (\cite{ceyhan:dom-num-NPE-SPL}),
the extension for general $M$ is similar:
Let $\Y = \{\y_1,\y_2,\ldots,\y_{d+1}\}$ be $d+1$ non-coplanar points.
Denote the simplex formed by these $d+1$ points as $\mathfrak S(\Y_{d+1})$.
The extension of $N^{\tau}_{CS}$ to $\R^d$ for $d > 2$ is straightforward.
Let $\Y = \{\y_1,\y_2,\cdots,\y_{d+1}\}$ be $d+1$
points in general position.
Denote the simplex formed by these $d+1$
points as $\mathfrak S(\Y_{d+1})$.
(A simplex is the simplest polytope in
$\R^d$ having $d+1$ vertices, $d\,(d+1)/2$ edges and $d+1$ faces of dimension $(d-1)$.)
For $\tau \in [0,1]$, define the central
similarity proximity map as follows. Let $\varphi_j$ be the face
opposite vertex $\y_j$ for $j=1,2,\ldots,d+1$, and ``face regions''
$R(\varphi_1),\ldots,R(\varphi_{d+1})$ partition $\mathfrak S(\Y_{d+1})$
into $d+1$ regions, namely the $d+1$ polytopes with vertices being
the center of mass together with $d$ vertices chosen from $d+1$ vertices.
For $x \in \mathfrak S(\Y_{d+1}) \setminus \Y$, let $\varphi(x)$
be the face in whose region $x$ falls; $x \in R(\varphi(x))$.
(If $x$ falls on the boundary of two face regions,
we assign $\varphi(x)$ arbitrarily.)
For ${\tau} \in (0,1]$, the {\em $\tau$-factor} central similarity
proximity region $N(x,\tau,M_C)=N^{\tau}(x)$ is defined to be the
simplex $\mathfrak S_{\tau}(x)$ with the following properties:
\begin{itemize}
\item[(i)] $\mathfrak S_{\tau}(x)$ has a face $\varphi_{\tau}(x)$ parallel to $\varphi(x)$ such that
$\tau\, d(x,\varphi(x))=d(\varphi_{\tau}(x),x)$ where $d(x,\varphi(x))$ is the Euclidean (perpendicular) distance from $x$ to $\varphi(x)$ ,
\item[(ii)] $\mathfrak S_{\tau}(x)$ has the same orientation as and is similar to $\mathfrak S(\Y_{d+1})$,
\item[(iii)] $x$ is at the center of mass of $\mathfrak S_{\tau}(x)$.  Note that $\tau>1$ implies that $x \in N(x,\tau,M_C)$.
\end{itemize}
For $\tau=0$, define $N(x,\tau,M_C)=\{x\}$ for all $x \in \mathfrak S(\Y_{d+1})$.

Theorem \ref{thm:scale-inv-NCSt} generalizes, so that any simplex $\mathfrak S$ in $\R^d$
can be transformed into a regular polytope (with edges being equal
in length and faces being equal in volume) preserving uniformity.
Delaunay triangulation becomes Delaunay tessellation in $\R^d$,
provided no more than $d+1$ points being cospherical (lying on the
boundary of the same sphere).
In particular, with $d=3$, the general
simplex is a tetrahedron (4 vertices, 4 triangular faces and 6
edges), which can be mapped into a regular tetrahedron (4 faces are
equilateral triangles) with vertices
$(0,0,0)\,(1,0,0)\,(1/2,\sqrt{3}/2,0),\,(1/2,\sqrt{3}/6,\sqrt{6}/3)$.

Asymptotic normality of the $U$-statistic and consistency of the
tests hold for $d>2$.

\section{Discussion}
\label{sec:disc-conclusions}
In this article,
we consider the relative density of a
random digraph family called central similarity proximity catch digraph (PCD)
which is based on two classes of points (in $\mathbb R$).
The central similarity PCDs have an expansion parameter $\tau >0$
and
a centrality parameter $c \in (0,1/2)$.
We demonstrate that the relative density
of the central similarity PCDs is a $U$-statistic.
Then, applying the central limit theory of the $U$-statistics,
we derive the (asymptotic normal) distribution of
the relative density for uniform data for the entire ranges of $\tau$ and $c$.
%Relative density is the number of arcs in a given digraph to the total number of arcs possible
%in a complete symmetric digraph with the same number of vertices.
%Points from one of the classes constitute
%the vertices of the PCDs and are a random sample from uniform distribution
%in compact intervals in $\R$.
%We provide the asymptotic distribution of
%the relative density for central similarity PCDs
%for uniform data .
We also determine the parameters $\tau$ and $c$ for which
the rate of convergence to normality is the fastest.
%The PCD in this article can also be viewed as the one dimensional version
%of the PCD in \cite{ceyhan:dom-num-NPE-SPL,ceyhan:dom-num-NPE-MASA}
%(see also Section \ref{sec:NCSt-higher-D}).
%These generalizations are parameterized with two parameters:
%an expansion parameter $r$ and a centrality parameter $c$.

We can apply the relative density in testing
one dimensional bivariate spatial point patterns,
as done in \cite{ceyhan:arc-density-CS} for two-dimensional data.
Let $\X$ and $\Y$ be two classes of points which lie in a compact interval in $\mathbb R$.
Then our null hypothesis is some form of complete spatial randomness
of $\X$ points, which implies that distribution of $\X$ points
has a uniform distribution in the support interval irrespective of the distribution of the $\Y$ points.
The alternatives are the segregation of $\X$ from $\Y$ points
or association of $\X$ points with $\Y$ points.
In general, association is the pattern in which
the points from the two different classes occur close to each other,
while
segregation is the pattern in which
the points from the same class tend to cluster together.
In this context,
under association, $\X$ points are clustered around $\Y$ points,
while
under segregation, $\X$ points are clustered away from the $\Y$ points.
Notice that we can use the asymptotic distribution (i.e., the normal approximation)
of the relative density for spatial pattern tests,
so our methodology requires number of $\X$ points to be much larger compared to the number of $\Y$ points.
Our results will make the power comparisons possible
for data from large families of distributions.
Moreover, one might determine the optimal (with respect to empirical size and power)
parameter values against segregation and association alternatives.

The central similarity PCDs for one dimensional data can be used in classification
as outlined in \cite{priebe:2003b},
if a high dimensional data set can be projected to one dimensional space with unsubstantial information loss
(by some dimension reduction method).
In the classification procedure, one might also determine
the optimal parameters (with respect to some penalty function) for the best performance.
Furthermore, this work forms the foundation of the generalizations and calculations
for uniform and non-uniform cases in multiple dimensions.
See Section \ref{sec:NCSt-higher-D}
for the details of the extension to higher dimensions.
For example, in $\mathbb R^2$,
the expansion parameter is still $\tau$,
but the centrality parameter is $M=(m_1,m_2)$,
which is two dimensional.
The optimal parameters for
testing spatial patterns and classification can also be determined,
as in the one dimensional case.

\section*{Acknowledgments}
%I would like to thank the anonymous referees whose constructive
%comments and suggestions greatly improved the presentation and flow
%of this article.
This work was supported by TUBITAK Kariyer Project Grant 107T647.

%\bibliography{References}
%\bibliographystyle{apalike}
%\bibliographystyle{plain}

\subsection*{APPENDIX 1: Proofs for the One Interval Case}

\subsection*{Proof of Theorem \ref{thm:reldens-tau=1,c}:}
Depending on the location of $x_1$,
the following are the different types of the combinations of $N(x_1,1,c)$ and $\G_1(x_1,1,c)$.
\begin{itemize}
\item[(i)] for $0< x_1 \le c$, we have $N(x_1,1,c)=(0,x_1/c)$ and $\G_1(x_1,1,c)=(c\,x_1,(1-c)\,x_1+c)$,
\item[(ii)] for $c < x_1 < 1$, $N(x_1,1,c)=((x_1-c)/(1-c),1)$ and $\G_1(x_1,1,c)=(c\,x_1,(1-c)\,x_1+c)$.
\end{itemize}

Then
$\mu(1,c)=P(X_2 \in N(X_1,1,c))
=\int_0^{c} \frac{x_1}{c} dx_1+\int_c^1 (1-\frac{x_1-c}{1-c}) dx_1
=1/2$.

For $\Cov(h_{12},h_{13})$,
we need to calculate $P_{2N}$, $P_{NG}$, and $P_{2G}$.
$$P_{2N}=P(\{X_2,X_3\} \subset N(X_1,1,c))
=\int_0^{c} \left(\frac{x_1}{c}\right)^2 dx_1+\int_c^1 \left(1-\frac{x_1-c}{1-c}\right)^2 dx_1
=1/3.$$

\begin{multline*}
P_{NG}=P(X_2 \in N(X_1,1,c),X_3 \in \G_1(X_1,1,c))=\\
\int_0^{c} \frac{x_1}{c}(1+c-2\,c\,x_1) dx_1+\int_c^1 \left(1-\frac{x_1-c}{1-c}\right)(1+c-2\,c\,x_1) dx_1
=-c^2/3+c/3+1/6.
\end{multline*}

Finally,
$P_{2G}=P(\{X_2,X_3\} \subset \G_1(X_1,1,c))
=\int_0^1 (1+c-2\,c\,x_1)^2 dx_1
=c^2/3-c/3+1/3.$

Therefore
$4\,\E[h_{12}h_{13}]=P_{2N}+2\,P_{NG}+P_{2G}=
-c^2/3+c/3+1$.
Hence
$4\,\nu(1,c)=4\,\Cov[h_{12},h_{13}]=c\,(1-c)/3$.
$\blacksquare$

\subsection*{Proof of Theorem \ref{thm:reldens-tau,c=1/2}:}
There are two cases for $\tau$,
namely $0 < \tau < 1$ and $\tau \ge 1$.\\
\textbf{Case 1:} $0 < \tau < 1$:
In this case depending on the location of $x_1$,
the following are the different types of the combinations of $N(x_1,\tau,1/2)$ and $\G_1(x_1,\tau,1/2)$.
\begin{itemize}
\item[(i)] for $0 < x_1 \le (1-\tau)/2$, we have
$N(x_1,\tau,1/2)=(x_1\,(1-\tau),x_1\,(1+\tau))$ and $\G_1(x_1,\tau,1/2)=(x_1/(1+\tau),x_1/(1-\tau))$,
\item[(ii)] for $(1-\tau)/2 < x_1 \le 1/2$, we have
$N(x_1,\tau,1/2)=(x_1\,(1-\tau),x_1\,(1+\tau))$ and $\G_1(x_1,\tau,1/2)=(x_1/(1+\tau),(x_1+\tau)/(1+\tau))$.
\end{itemize}

Then
$\mu(\tau,1/2)=P(X_2 \in N(X_1,\tau,1/2))=2\,P(X_2 \in N(X_1,\tau,1/2),X_1 \in (0,1/2))$ by symmetry
and
$$P(X_2 \in N(X_1,\tau,1/2),X_1 \in (0,1/2))=
\int_0^{1/2} (x_1\,(1+\tau)-x_1\,(1-\tau)) dx_1=
\int_0^{1/2} 2\,x_1\,\tau dx_1=
\tau/4.$$
So
$\mu(\tau,1/2)=2\,(\tau/4)=\tau/2$.

For $\Cov(h_{12},h_{13})$, we need to calculate $P_{2N}$, $P_{NG}$, and $P_{2G}$.
$$P_{2N}=P(\{X_2,X_3\} \subset N(X_1,\tau,1/2))=2\,P(\{X_2,X_3\} \subset N(X_1,\tau,1/2),X_1 \in (0,1/2))$$
and
$$P(\{X_2,X_3\} \subset N(X_1,\tau,1/2),X_1 \in (0,1/2))
=\int_0^{1/2} (2\,x_1\,\tau)^2 dx_1=\tau^2/6.$$
So
$P_{2N}=2\,(\tau^2/6)=\tau^2/3$.

\begin{multline*}
P_{NG}=P(X_2 \in N(X_1,\tau,1/2),X_3 \in \G_1(X_1,\tau,1/2))=\\
2\,P(X_2 \in N(X_1,\tau,1/2),X_3 \in \G_1(X_1,\tau,1/2),X_1 \in (0,1/2))
\end{multline*}
and
\begin{multline*}
P(X_2 \in N(X_1,\tau,1/2),X_3 \in \G_1(X_1,\tau,1/2),X_1 \in (0,1/2)=\\
\int_0^{(1-\tau)/2} (2\,x_1\,\tau)\left(\frac{x_1}{1-\tau}-\frac{x_1}{1+\tau}\right) dx_1+
\int_{(1-\tau)/2}^{1/2}(2\,x_1\,\tau)\left(\frac{x_1+\tau}{1+\tau}-\frac{x_1}{1+\tau}\right) dx_1=\\
\int_0^{(1-\tau)/2} (2\,x_1\,\tau)\left(\frac{2\,x_1\,\tau}{1-\tau^2}\right) dx_1+
\int_{(1-\tau)/2}^{1/2}(2\,x_1\,\tau)\left(\frac{\tau}{1+\tau}\right) dx_1
={\frac { \left( 2+2\,\tau-\tau^2 \right) \tau^2}{12\,(\tau+1)}}.
\end{multline*}
So
$P_{NG}={\frac { \left( 2+2\,\tau-\tau^2  \right) \tau^2}{6\,(\tau+1)}}$.

Finally,
$$P_{2G}=P(\{X_2,X_3\} \subset \G_1(X_1,\tau,1/2))=2\,P(\{X_2,X_3\} \subset \G_1(X_1,\tau,1/2),X_1 \in (0,1/2))$$
and
$$P(\{X_2,X_3\} \subset \G_1(X_1,\tau,1/2),X_1 \in (0,1/2))
=\int_0^{(1-\tau)/2} \left(\frac{2\,x_1\,\tau}{1-\tau^2}\right)^2 dx_1+
\int_{(1-\tau)/2}^{1/2}\left(\frac{\tau}{1+\tau}\right)^2 dx_1
={\frac {\tau^2 \left( 2\,\tau+1 \right) }{6\,\left( \tau+1 \right)^2}}.$$
So
$P_{2G}={\frac {\tau^2 \left( 2\,\tau+1 \right) }{ 6\,\left( \tau+1 \right)^2}}$.

Therefore
$4\,\E[h_{12}h_{13}]=P_{2N}+2\,P_{NG}+P_{2G}=
{\frac {\tau^2 \left( 8\,\tau+4-\tau^3+2\,\tau^2 \right) }{ 3\,\left( \tau+1 \right)^2}}$.
Hence
$4\,\nu(\tau,1/2)=4\,\Cov[h_{12},h_{13}]={\frac {\tau^2 \left( -\tau^3-\tau^2+2\,\tau+1 \right) }{ 3\,\left( \tau+1 \right)^2}}$.

\textbf{Case 2:} $\tau \ge 1$:
In this case depending on the location of $x_1$,
the following are the different types of the combinations of $N(x_1,\tau,1/2)$ and $\G_1(x_1,\tau,1/2)$.
\begin{itemize}
\item[(i)] for $0 < x_1 \le 1/(1+\tau)$,  we have
$N(x_1,\tau,1/2)=(0,x_1\,(1+\tau))$ and $\G_1(x_1,\tau,1/2)=(x_1/(1+\tau),(x_1+\tau)/(1+\tau))$,
\item[(ii)] for $1/(1+\tau) < x_1 \le 1/2$,  we have
$N(x_1,\tau,1/2)=(0,1)$ and $\G_1(x_1,\tau,1/2)=(x_1/(1+\tau),(x_1+\tau)/(1+\tau))$,
\end{itemize}

Then
$\mu(\tau,1/2)=P(X_2 \in N(X_1,\tau,1/2))=2\,P(X_2 \in N(X_1,\tau,1/2),X_1 \in (0,1/2))$ by symmetry
and
$$P(X_2 \in N(X_1,\tau,1/2),X_1 \in (0,1/2))
=\int_0^{1/(1+\tau)} x_1\,(1+\tau) dx_1+\int_{1/(1+\tau)}^{1/2} 1 dx_1=\frac{\tau}{2\,(\tau+1)}.$$
So
$\mu(\tau,1/2)=2\,\left(\frac{\tau}{2\,(\tau+1)}\right)=\frac{\tau}{(\tau+1)}$.

Next
$$P_{2N}=P(\{X_2,X_3\} \subset N(X_1,\tau,1/2))=2\,P(\{X_2,X_3\} \subset N(X_1,\tau,1/2),X_1 \in (0,1/2))$$
and
$$P(\{X_2,X_3\} \subset N(X_1,\tau,1/2),X_1 \in (0,1/2))
=\int_0^{1/(1+\tau)} (x_1\,(1+\tau))^2 dx_1+\int_{1/(1+\tau)}^{1/2} 1 dx_1
={\frac {1-3\,\tau}{6\,(\tau+1)}}.$$
So
$P_{2N}=2\,\left({\frac {1-3\,\tau}{6\,(\tau+1)}}\right)={\frac {1-3\,\tau}{3\,(\tau+1)}}$.

\begin{multline*}
P_{NG}=P(X_2 \in N(X_1,\tau,1/2),X_3 \in \G_1(X_1,\tau,1/2))=\\
2\,P(X_2 \in N(X_1,\tau,1/2),X_3 \in \G_1(X_1,\tau,1/2),X_1 \in (0,1/2))
\end{multline*}
and
\begin{multline*}
P(X_2 \in N(X_1,\tau,1/2),X_3 \in \G_1(X_1,\tau,1/2),X_1 \in (0,1/2)=\\
\int_0^{1/(1+\tau)} (x_1\,(1+\tau))(\tau/(1+\tau)) dx_1+\int_{1/(1+\tau)}^{1/2} (\tau/(1+\tau)) dx_1
=\frac{\tau^2}{2\,(1+\tau)^2}.
\end{multline*}
So
$P_{NG}=\frac{\tau^2}{(1+\tau)^2}$.

Finally,
$$P_{2G}=P(\{X_2,X_3\} \subset \G_1(X_1,\tau,1/2))=2\,P(\{X_2,X_3\} \subset \G_1(X_1,\tau,1/2),X_1 \in (0,1/2))$$
and
$$P(\{X_2,X_3\} \subset \G_1(X_1,\tau,1/2),X_1 \in (0,1/2))
=\int_0^{1/2} (\tau/(1+\tau))^2 dx_1
=\frac{\tau^2}{2\,(1+\tau)^2}.$$
So
$P_{2G}=\frac{\tau^2}{(1+\tau)^2}$.

Therefore
$4\,\E[h_{12}h_{13}]=P_{2N}+2\,P_{NG}+P_{2G}={\frac {12\,\tau^2+2\,\tau-1}{ 3\,\left( \tau+1 \right)^2}}$.
Hence
$4\,\nu(\tau,1/2)=4\,\Cov[h_{12},h_{13}]={\frac {2\,\tau-1}{3\,\left( \tau+1 \right)^2}}$.
$\blacksquare$

\subsection*{Proof of Theorem \ref{thm:reldens-tau,c}:}
First we consider $0 < c \le 1/2$.
There are two cases for $\tau$,
namely $0 < \tau < 1$ and $ \tau \ge 1$.\\
\textbf{Case 1:} $0 < \tau < 1$:
In this case depending on the location of $x_1$,
the following are the different types of the combinations of $N(x_1,\tau,c)$ and $\G_1(x_1,\tau,c)$.
Let
$a_1:=x_1\,(1-\tau)$,
$a_2:=x_1\,(1+\frac{(1-c)\,\tau}{c})$,
$a_3:=x_1-\frac{c\,\tau\,(1-x_1)}{1-c}$,
$a_4:=x_1+(1-x_1)\,\tau$,
and
$g_1:=\frac{c\,x_1}{c+(1-c)\,\tau}$,
$g_2:=\frac{x_1}{1-\tau}$,
$g_3:=\frac{x_1-\tau}{1-\tau}$,
$g_4:=\frac{x_1\,(1-c)+c\,\tau}{1-c+c\,\tau}$.
Then
\begin{itemize}
\item[(i)] for $0 < x_1 \le c\,(1-\tau)$, we have
$N(x_1,\tau,c)=(a_1,a_2)$ and $\G_1(x_1,\tau,c)=(g_1,g_2)$,
\item[(ii)] for $c\,(1-\tau) < x_1 \le c$, we have
$N(x_1,\tau,c)=(a_1,a_2)$ and $\G_1(x_1,\tau,c)=(g_1,g_4)$,
\item[(iii)] for $c < x_1 \le c\,(1-\tau)+\tau$, we have
$N(x_1,\tau,c)=(a_3,a_4)$ and $\G_1(x_1,\tau,c)=(g_1,g_4)$,
\item[(iv)] for $c\,(1-\tau)+\tau < x_1 < 1$, we have
$N(x_1,\tau,c)=(a_3,a_4)$ and $\G_1(x_1,\tau,c)=(g_3,g_4)$.
\end{itemize}

Then
$\mu(\tau,c)=P(X_2 \in N(X_1,\tau,c))=
\int_0^{c} (a_2-a_1) dx_1+\int_c^1 (a_4-a_3) dx_1
=\tau/2$.

For $\Cov(h_{12},h_{13})$, we need to calculate $P_{2N}$, $P_{NG}$, and $P_{2G}$.
$$P_{2N}=P(\{X_2,X_3\} \subset N(X_1,\tau,c))
=\int_0^{c} (a_2-a_1)^2 dx_1+\int_c^1 (a_4-a_3)^2 dx_1
=\tau^2/3.$$

\begin{multline*}
P_{NG}=P(X_2 \in N(X_1,\tau,c),X_3 \in \G_1(X_1,\tau,c))=
\int_0^{c\,(1-\tau)} (a_2-a_1)(g_2-g_1) dx_1+\\
\int_{c\,(1-\tau)}^{c} (a_2-a_1)(g_4-g_1) dx_1+
\int_c^{c\,(1-\tau)+\tau} (a_4-a_3)(g_4-g_1) dx_1+
\int_{c\,(1-\tau)+\tau}^1 (a_4-a_3)(g_4-g_3) dx_1=\\
{\frac {\tau^2 \left( c^2\,\tau^3-5\,c^2\,\tau^2
-c\,\tau^3+4\,c^2\,\tau+5\,c\,\tau^2-2\,c^2-4\,c\,\tau-\tau^2+2\,c+2\,\tau \right) }
{6\,\left( c\,\tau-c+1 \right)  \left(c+\tau-c\,\tau \right) }}.
\end{multline*}
Finally,
\begin{multline*}
P_{2G}=P(\{X_2,X_3\} \subset \G_1(X_1,\tau,c))=
\int_0^{c\,(1-\tau)} (g_2-g_1)^2 dx_1+\\
\int_{c\,(1-\tau)}^{c\,(1-\tau)+\tau} (g_4-g_1)^2 dx_1+
\int_{c\,(1-\tau)+\tau}^1 (g_4-g_3)^2 dx_1
={\frac { \left( 2\,c^2\,\tau-c^2-2\,c\,\tau+c+\tau \right) \tau^2}{ 3\,\left( c\,\tau-c+1 \right)  \left( c+\tau-c\,\tau \right) }}.
\end{multline*}

Therefore
$$4\,\E[h_{12}h_{13}]=P_{2N}+2\,P_{NG}+P_{2G}=
{\frac {\tau^2 \left( c^2\,\tau^3-6\,c^2\,\tau^2-
c\,\tau^3+8\,c^2\,\tau+6\,c\,\tau^2-4\,c^2-8\,c\,\tau-\tau^2+
4\,c+4\,\tau \right) }{3\,\left( c\,\tau-c+1 \right)  \left( c+\tau-c\,\tau \right) }}.$$
Hence
$4\,\kappa_1(\tau,c)=4\,\Cov[h_{12},h_{13}]=
{\frac {\tau^2 \left( c^2\,\tau^3-3\,c^2\,\tau^2-
c\,\tau^3+2\,c^2\,\tau+3\,c\,\tau^2-c^2-2\,c\,\tau-\tau^2+
c+\tau \right) }{ 3\,\left( c\,\tau-c+1 \right)  \left( c+\tau-c\,\tau
 \right) }}$.

\textbf{Case 2:} $\tau \ge 1$:
In this case depending on the location of $x_1$,
the following are the different types of the combinations of $N(x_1,\tau,c)$ and $\G_1(x_1,\tau,c)$.
\begin{itemize}
\item[(i)] for $0 < x_1 \le \frac{c}{c+(1-c)\,\tau}$, we have
$N(x_1,\tau,c)=(0,a_2)$ and $\G_1(x_1,\tau,c)=(g_1,g_4)$,
\item[(ii)] for $\frac{c}{c+(1-c)\,\tau} < x_1 \le  \frac{c\,\tau}{1-c+c\,\tau}$, we have
$N(x_1,\tau,c)=(0,1)$ and $\G_1(x_1,\tau,c)=(g_1,g_4)$,
\item[(iii)] for $\frac{c\,\tau}{1-c+c\,\tau} < x_1 < 1$, we have
$N(x_1,\tau,c)=(a_3,1)$ and $\G_1(x_1,\tau,c)=(g_1,g_4)$.
\end{itemize}

Then
\begin{multline*}
\mu(\tau,c)=P(X_2 \in N(X_1,\tau,c))=
\int_0^{\frac{c}{c+(1-c)\,\tau}} a_2 dx_1+
\int_{\frac{c}{c+(1-c)\,\tau}}^{\frac{c\,\tau}{1-c+c\,\tau}} 1 dx_1+
\int_{\frac{c\,\tau}{1-c+c\,\tau}}^1 (1-a_3) dx_1=\\
{\frac {\tau\, \left( 2\,c^2\,\tau-2\,c^2-2\,c\,\tau+2\,c-1\right) }
{2\,\left( c\,\tau-c+1 \right)  \left( c\,\tau-c-\tau \right) }}.
\end{multline*}

Next
\begin{multline*}
P_{2N}=P(\{X_2,X_3\} \subset N(X_1,\tau,c))=
\int_0^{\frac{c}{c+(1-c)\,\tau}} a_2^2 dx_1+
\int_{\frac{c}{c+(1-c)\,\tau}}^{\frac{c\,\tau}{1-c+c\,\tau}} 1 dx_1+
\int_{\frac{c\,\tau}{1-c+c\,\tau}}^1 (1-a_3)^2 dx_1=\\
{\frac {3\,c^2\,\tau^2-2\,c^2\,\tau-3\,c\,\tau^2-c^2+
2\,c\,\tau+c-\tau}{ 3\,\left( c\,\tau-c+1 \right)  \left( c\,\tau-c-\tau \right) }}.
\end{multline*}

\begin{multline*}
P_{NG}=P(X_2 \in N(X_1,\tau,c),X_3 \in \G_1(X_1,\tau,c))=\\
\int_0^{\frac{c}{c+(1-c)\,\tau}} a_2\,(g_4-g_1) dx_1+
\int_{\frac{c}{c+(1-c)\,\tau}}^{\frac{c\,\tau}{1-c+c\,\tau}} (g_4-g_1) dx_1+
\int_{\frac{c\,\tau}{1-c+c\,\tau}}^1 (1-a_3)\,(g_4-g_1) dx_1=\\
[\tau^2 \bigl( 6\,c^6\,\tau^4-24\,c^6\,\tau^3-18\,c^5\,\tau^4+36\,c^6\,\tau^2+72\,c^5\,\tau^3+
18\,c^4\,\tau^4-24\,c^6\,\tau-108\,c^5\,\tau^2-84\,c^4\,\tau^3-
6\,c^3\,\tau^4+6\,c^6+\\
72\,c^5\,\tau+132\,c^4\,\tau^2+
48\,c^3\,\tau^3-18\,c^5-92\,c^4\,\tau-84\,c^3\,\tau^2-
12\,c^2\,\tau^3+26\,c^4+64\,c^3\,\tau+30\,
c^2\,\tau^2-22\,c^3-26\,c^2\,\tau-6\,c\,\tau^2+\\
10\,c^2+6\,c\,\tau-2\,c-\tau \bigr) ]/[6\,\left( c\,\tau-c+1 \right)^3 \left(c\,\tau-c-\tau \right)^3].
\end{multline*}

Finally,
\begin{multline*}
P_{2G}=P(\{X_2,X_3\} \subset \G_1(X_1,\tau,c))
=\int_0^1 (g_4-g_1)^2 dx_1=\\
{\frac { \left( 3\,c^4\,\tau^2-6\,c^4\,\tau-6\,c^3\,\tau^2+
3\,c^4+12\,c^3\,\tau+3\,c^2\,\tau^2-6\,c^3-9
\,c^2\,\tau+7\,c^2+3\,c\,\tau-4\,c+1 \right) \tau^2}
{3\, \left( c\,\tau-c+1 \right)^2 \left( c\,\tau-c-\tau \right)^2}}.
\end{multline*}

Therefore
\begin{multline*}
4\,\E[h_{12}h_{13}]=P_{2N}+2\,P_{NG}+P_{2G}=
[12\,c^6\,\tau^6-50\,c^6\,\tau^5-36\,c^5\,\tau^6+
79\,c^6\,\tau^4+150\,c^5\,\tau^5+36\,c^4\,\tau^6-56\,c^6\,\tau^3-\\
237\,c^5\,\tau^4-175\,c^4\,\tau^5-
12\,c^3\,\tau^6+14\,c^6\,\tau^2+168\,c^5\,\tau^3+
297\,c^4\,\tau^4+100\,c^3\,\tau^5+2\,c^6\,\tau-
42\,c^5\,\tau^2-220\,c^4\,\tau^3-199\,c^3\,\tau^4-\\
25\,c^2\,\tau^5-c^6-6\,c^5\,\tau+58\,c^4\,\tau^2+
160\,c^3\,\tau^3+75\,c^2\,\tau^4+3\,c^5+7\,c^4\,\tau-
46\,c^3\,\tau^2-70\,c^2\,\tau^3-15\,c\,\tau^4-3\,c^4-
4\,c^3\,\tau+20\,c^2\,\tau^2+\\
18\,c\,\tau^3+c^3+c^2\,\tau-4\,c\,\tau^2-3\,\tau^3]/
[ 3\,\left( c\,\tau-c+1 \right)^3\left( c\,\tau-c-\tau \right)^3].
\end{multline*}
Hence
\begin{multline*}
4\,\kappa_2(\tau,c)=4\,\Cov[h_{12},h_{13}]=
[c \left( 1-c \right)  \bigl( 2\,c^4\,\tau^5-7\,c^4\,\tau^4-
4\,c^3\,\tau^5+8\,c^4\,\tau^3+14\,c^3\,\tau^4+
3\,c^2\,\tau^5-2\,c^4\,\tau^2-\\
16\,c^3\,\tau^3-7\,c^2\,\tau^4-c\,\tau^5-2\,c^4\,\tau+4\,c^3\,\tau^2+
12\,c^2\,\tau^3+c^4+4\,c^3\,\tau-6\,c^2\,\tau^2-
4\,c\,\tau^3-2\,c^3-3\,c^2\,\tau+4\,c\,\tau^2+c^2+c\,\tau-\tau^2 \bigr) ]/\\
[ 3\,\left( c\,\tau-c+1 \right)^3 \left( c\,\tau-c-\tau \right)^3].
\end{multline*}

For $1/2 \le c < 1$,
by symmetry, it follows that
$\mu_2(\tau,c)=\mu_1(\tau,1-c)$
and
$\nu_2(\tau,c)=\nu_1(\tau,1-c)$.
$\blacksquare$

\subsection*{Proof of Theorem \ref{thm:reldens-end-int-tau=1}:}
Suppose $i=m+1$ (i.e., the support is the right end interval).
For $x_1 \in (0,1)$,
depending on the location of $x_1$,
the following are the different types of the combinations of $N_e(x_1,1)$ and $\G_{1,e}(x_1,1)$.
\begin{itemize}
\item[(i)] for $0<x_1 \le 1/2$,  we have
$N_e(x_1,1)=(0,2\,x_1)$ and $\G_{1,e}(x_1,1)=(x_1/2,1)$,
\item[(ii)] for $1/2<x_1<1$,
$N_e(x_1,1)=(0,1)$ and $\G_{1,e}(x_1,1)=(x_1/2,1)$.
\end{itemize}
Then
$\mu_e(1)=P(X_2 \in N_e(X_1,1))
=\int_0^{1/2} 2x_1 dx_1+\int_{1/2}^1 1 dx_1 =3/4$.

For $\Cov(h_{12},h_{13})$, we need to calculate $P_{2N}$, $P_{NG}$, and $P_{2G}$.
$$P_{2N}=P(\{X_2,X_3\} \subset N_e(X_1,1))
=\int_0^{1/2} (2x_1)^2 dx_1+\int_{1/2}^1 1 dx_1 =2/3.$$

$$P_{NG}=P(X_2 \in N_e(X_1,1),X_3 \in \G_{1,e}(X_1,1))
=\int_0^{1/2} (2x_1)(1-x_1/2) dx_1+\int_{1/2}^{1} 1(1-x_1/2) dx_1
=25/48.$$

Finally,
$P_{2G}=P(\{X_2,X_3\} \subset \G_{1,e}(X_1,1))
=\int_0^{1} (1-x_1/2)^2 dx_1=7/12.$

Therefore
$4\,\E[h_{12}h_{13}]=P_{2N}+2\,P_{NG}+P_{2G}=55/24$.
Hence
$4\,\nu_e(1)=4\,\Cov[h_{12},h_{13}]=1/24$.

For uniform data,
by symmetry,
the distribution of the relative density of the subdigraph for $i=1$ is identical to $i=m+1$ case.
$\blacksquare$

\subsection*{Proof of Theorem \ref{thm:reldens-end-int-tau}:}
There are two cases for $\tau$,
namely, $0<\tau<1$ and $\tau \ge 1$.

\textbf{Case 1:} $0 < \tau < 1$:
For $x_1 \in (0,1)$,
depending on the location of $x_1$,
the following are the different types of the combinations of $N_e(x_1,\tau)$ and $\G_{1,e}(x_1,\tau)$.
\begin{itemize}
\item[(i)] for $0 < x_1 \le 1-\tau$, we have
$N_e(x_1,\tau)=(x_1\,(1-\tau),x_1\,(1+\tau))$ and $\G_{1,e}(x_1,\tau)=(x_1/(1+\tau),x_1/(1-\tau))$,
\item[(ii)] for $1-\tau < x_1 \le 1/(1+\tau)$, we have
$N_e(x_1,\tau)=(x_1\,(1-\tau),x_1\,(1+\tau))$ and $\G_{1,e}(x_1,\tau)=(x_1/(1+\tau),1)$,
\item[(iii)] for $1/(1+\tau) < x_1 < 1$, we have
$N_e(x_1,\tau)=(x_1\,(1-\tau),1)$ and $\G_{1,e}(x_1,\tau)=(x_1/(1+\tau),1)$.
\end{itemize}
Then
\begin{multline*}
\mu_e(\tau)=P(X_2 \in N_e(X_1,\tau))=
\int_0^{1/(1+\tau)} (x_1\,(1+\tau)-x_1\,(1-\tau)) dx_1+
\int_{1/(1+\tau)}^1 (1-x_1\,(1-\tau)) dx_1 =\\
\int_0^{1/(1+\tau)} (2\,x_1\,\tau) dx_1+
\int_{1/(1+\tau)}^1 (1-x_1+x_1\,\tau) dx_1 =
\frac{\tau\,(\tau+2)}{2\,(\tau+1)}.
\end{multline*}

For $\Cov(h_{12},h_{13})$, we need to calculate $P_{2N,e}$, $P_{NG,e}$, and $P_{2G,e}$.
$$P_{2N,e}=P(\{X_2,X_3\} \subset N_e(X_1,\tau))
=\int_0^{1/(1+\tau)} (2\,x_1\,\tau)^2 dx_1+\int_{1/(1+\tau)}^1 (1-x_1+x_1\,\tau)^2 dx_1 =
\frac{\tau^2\,(\tau^2+3\,\tau+4)}{3\,(\tau+1)^2}.$$

\begin{multline*}
P_{NG,e}=P(X_2 \in N_e(X_1,\tau),X_3 \in \G_{1,e}(X_1,\tau))=
\int_0^{1-\tau} (2\,x_1\,\tau)\left(\frac{2\,x_1\,\tau}{1-\tau^2}\right) dx_1+\\
\int_{1-\tau}^{1/(1+\tau)} (2\,x_1\,\tau)\left(1-\frac{x_1}{1+\tau}\right) dx_1+
\int_{1/(1+\tau)}^1 (1-x_1\,(1-\tau))\left(1-\frac{x_1}{1+\tau}\right) dx_1
={\frac { \left( 7\,{\tau}^2+14\,\tau+8-2\,{\tau}^4-2\,{\tau}^3 \right) {\tau}^2}
{ 6\,\left( \tau+1 \right)^3}}.
\end{multline*}

Finally,
$$P_{2G,e}=P(\{X_2,X_3\} \subset \G_{1,e}(X_1,\tau))
=\int_0^{1-\tau} \left(\frac{2\,x_1\,\tau}{1-\tau^2}\right)^2 dx_1+
\int_{1-\tau}^1\left(1-\frac{x_1}{1+\tau}\right)^2 dx_1=
\frac{\tau^2\,(3\,\tau+4)}{3\,(\tau+1)^2}.$$

Therefore
$4\,\E[h_{12}h_{13}]=P_{2N,e}+2\,P_{NG,e}+P_{2G,e}=
{\frac {{\tau}^2 \left( 2\,{\tau}^2+5\,\tau+4 \right)
\left( 2\,\tau+4-{\tau}^2 \right) }{ 3\,\left( \tau+1 \right)^3}}$.
Hence
$$4\,\nu_e(\tau)=4\,\Cov[h_{12},h_{13}]=
{\frac {{\tau}^2 \left( 4\,\tau+4-2\,{\tau}^4-4\,{\tau}^3-{\tau}^2 \right) }
{ 3\,\left( \tau+1 \right)^3}}.$$

\textbf{Case 2:} $\tau \ge 1$:
For $x_1 \in (0,1)$,
depending on the location of $x_1$,
the following are the different types of the combinations of $N_e(x_1,\tau)$ and $\G_{1,e}(x_1,\tau)$.
\begin{itemize}
\item[(i)] for $0 < x_1 \le 1/(1+\tau)$, we have
$N_e(x_1,\tau)=(0,x_1\,(1+\tau))$ and $\G_{1,e}(x_1,\tau)=(x_1/(1+\tau),1)$,
\item[(ii)] for $1/(1+\tau) < x_1 < 1$, we have
$N_e(x_1,\tau)=(0,1)$ and $\G_{1,e}(x_1,\tau)=(x_1/(1+\tau),1)$.
\end{itemize}
Then
$$\mu_e(\tau)=P(X_2 \in N_e(X_1,\tau))
=\int_0^{1/(1+\tau)} x_1\,(1+\tau) dx_1+\int_{1/(1+\tau)}^1 1 dx_1 =
\frac{1+2\,\tau}{2\,(\tau+1)}.$$

Next,
$$P_{2N,e}=P(\{X_2,X_3\} \subset N_e(X_1,\tau))
=\int_0^{1/(1+\tau)} (x_1\,(1+\tau))^2 dx_1+\int_{1/(1+\tau)}^1 1 dx_1 =
\frac{1+3\,\tau}{3\,(\tau+1)}.$$

\begin{multline*}
P_{NG,e}=P(X_2 \in N_e(X_1,\tau),X_3 \in \G_{1,e}(X_1,\tau))=\\
\int_0^{1/(1+\tau)} (x_1\,(1+\tau))\left(1-\frac{x_1}{1+\tau}\right) dx_1+
\int_{1/(1+\tau)}^1 \left(1-\frac{x_1}{1+\tau}\right) dx_1
={\frac {6\,{\tau}^3+12\,{\tau}^2+6\,\tau+1}{ 6\,\left( \tau+1\right)^3}}.
\end{multline*}

Finally,
$$P_{2G,e}=P(\{X_2,X_3\} \subset \G_{1,e}(X_1,\tau))
=\int_0^1 \left(1-\frac{x_1}{1+\tau}\right)^2 dx_1=
\frac{3\,\tau^2+3\,\tau+1}{3\,(\tau+1)^2}.$$

Therefore
$4\,\E[h_{12}h_{13}]=P_{2N,e}+2\,P_{NG,e}+P_{2G,e}=
{\frac {12\,{\tau}^3+25\,{\tau}^2+15\,\tau+3}{ 3\,\left( \tau+1 \right)^3}}$.
Hence
$4\,\nu_e(\tau)=4\,\Cov[h_{12},h_{13}]=\frac{\tau^2}{3\,(\tau+1)^3}$.
$\blacksquare$

\subsection*{APPENDIX 2: Proofs for the Multiple Interval Case}
We give the proof of Theorem \ref{thm:MI-asy-norm-II} first.

\subsection*{Proof of Theorem \ref{thm:MI-asy-norm-II}:}
Recall that
$\widetilde{\rho}_{n,m}(\tau,c)$ is
the relative arc density of the PCD for the $m>2$ case.
Then it follows that
$\widetilde{\rho}_{n,m}(\tau,c)$ is a $U$-statistic of degree two,
so we can write it as
$\widetilde{\rho}_{n,m}(\tau,c)=\frac{2}{n(n-1)}\sum_{i<j}h_{ij}$
where
$h_{ij}=(g_{ij}+g_{ji})/2$.
Then the expectation of $\widetilde{\rho}_{n,m}(\tau,c)$ is
\begin{multline*}
\E\left[\widetilde{\rho}_{n,m}(\tau,c)\right]=
\frac{2}{n\,(n-1)}\sum\hspace*{-0.1 in}\sum_{i < j \hspace*{0.25 in}}
\hspace*{-0.1 in}\,\E\left[h_{ij}\right]=
\E\left[ h_{12} \right]=
\E\left[ g_{12} \right]=
P((X_1,X_2)\in \A)=
\widetilde \mu(m,\tau,c).
\end{multline*}
But, by definition of $N(\cdot,\tau,c)$,
if $X_1$ and $X_2$ are in different intervals,
then $P((X_1,X_2)\in \A)=0$.
So, by the law of total probability,
we have
\begin{multline*}
\widetilde \mu(m,\tau,c):=P((X_1,X_2)\in \A)=\\
\sum_{i=1}^{m+1} P((X_1,X_2)\in \A\,|\,\{X_1,X_2\} \subset \mI_i)\,P(\{X_1,X_2\} \subset \mI_i)=\\
\sum_{i=2}^{m} \mu(\tau,c)\,P(\{X_1,X_2\} \subset \mI_i)+\sum_{i \in \{1,m+1\}} \mu_e(\tau)\,P(\{X_1,X_2\} \subset \mI_i) = \\
\sum_{i=2}^{m} \mu(\tau,c)\,w_i^2+\sum_{i \in \{1,m+1\}} \mu_e(\tau)\,w_i^2
=\mu(\tau,c)\, \sum_{i=2}^{m} w_i^2+\mu_e(\tau)\,\sum_{i \in \{1,m+1\}} w_i^2.
\end{multline*}
since $P(X_2 \in N(X_1,\tau,c)\,|\,\{X_1,X_2\} \subset \mI_i)$ is $\mu(\tau,c)$ for middle intervals
and $\mu_e(\tau)$ for the end intervals
and $P(\{X_1,X_2\} \subset \mI_i)=\left(\frac{y_{(i)}-y_{(i-1)}}{\delta_2-\delta_1}\right)^2=w_i^2$.

Furthermore, the asymptotic variance is
$$4\,\widetilde \nu(m,\tau,c)=4\,\E\left[ h_{12}h_{13} \right]-\E\left[ h_{12} \right]\E\left[ h_{13} \right]=
4\,\E\left[ h_{12}h_{13} \right]-(\widetilde \mu(m,\tau,c))^2$$ where
$4\,\E\left[ h_{12}h_{13} \right]=\widetilde P_{2N}+2 \widetilde
P_{NG}+\widetilde P_{2G}$ with
\begin{multline*}
\widetilde P_{2N}=
\sum_{i=2}^{m} P(\{X_2,X_3\} \subset N(X_1,\tau,c)\,|\,\{X_1,X_2,X_3\} \subset \mI_i)\,P(\{X_1,X_2,X_3\} \subset \mI_i)+\\
\sum_{i\in \{1,m+1\}} P(\{X_2,X_3\} \subset N_e(X_1,\tau)\,|\,\{X_1,X_2,X_3\} \subset \mI_i)\,P(\{X_1,X_2,X_3\} \subset \mI_i)=\\
\sum_{i=2}^{m} P_{2N}\,P(\{X_1,X_2,X_3\} \subset \mI_i)+\sum_{i \in \{1,m+1\}}  P_{2N,e}\,P(\{X_1,X_2,X_3\} \subset \mI_i)\approx\\
\sum_{i=2}^{m} P_{2N}\,w_i^3+\sum_{i \in \{1,m+1\}} P_{2N,e}\,w_i^3
=P_{2N}\, \sum_{i=2}^{m} w_i^3+P_{2N,e}\,\sum_{i \in \{1,m+1\}} w_i^3.
\end{multline*}
since $P(\{X_2,X_3\} \subset N(X_1,\tau,c)\,|\,\{X_1,X_2,X_3\} \subset \mI_i)$ is $P_{2N}$ for middle intervals
and $P_{2N,e}$ for the end intervals
and $P(\{X_1,X_2,X_3\} \subset \mI_i)=\left(\frac{y_{(i)}-y_{(i-1)}}{\delta_2-\delta_1}\right)^3=w_i^3$.
Similarly,
$$
\widetilde P_{NG}=
P_{NG}\, \sum_{i=2}^{m} w_i^3+P_{NG,e}\,\sum_{i \in \{1,m+1\}} w_i^3
$$
and
$$
\widetilde P_{2G}=
P_{2G}\, \sum_{i=2}^{m} w_i^3+P_{2G,e}\,\sum_{i \in \{1,m+1\}} w_i^3.
$$
Therefore,
\begin{multline*}
4\,\widetilde \nu(m,\tau,c)=(P_{2N}+2\,P_{NG}+P_{2G})\sum_{i=2}^{m} w_i^3+
(P_{2N,e}+2\,P_{NG,e}+P_{2G,e})\sum_{i \in \{1,m+1\}} w_i^3
-(\widetilde \mu(m,\tau,c))^2.
\end{multline*}
Hence the desired result follows.
$\blacksquare$

\subsection*{Proof of Theorem \ref{thm:MI-asy-norm-I}:}
Recall that
$\rho_{n,m}(\tau,c)$ is the version I of the
relative arc density of the PCD for the $m>2$ case.
Moreover,
$\rho_{n,m}(\tau,c)=\frac{n(n-1)}{n_{_T}} \widetilde \rho_{n,m}(\tau,c)$.
Then the expectation of $\rho_{n,m}(\tau,c)$, for large $n_i$ and $n$, is

$$ \E\left[\rho_{n,m}(\tau,c)\right]=
\frac{n(n-1)}{n_{_T}} \E[\widetilde \rho_{n,m}(\tau,c)]
\approx \widetilde \mu(m,\tau,c) \, \left(\sum_{i=1}^{m+1} w_i^2\right)^{-1}$$
since
$\frac{n(n-1)}{n_{_T}}=\left(\sum_{i=1}^{m+1} n_i(n_i-1)/(n(n-1))\right)^{-1} \approx \left(\sum_{i=1}^{m+1} w_i^2\right)^{-1}$
for large $n_i$ and $n$.
Here $\widetilde \mu(m,\tau,c)$ is as in Theorem \ref{thm:MI-asy-norm-II}.

Moreover,
the asymptotic variance of $\rho_{n,m}(\tau,c)$, for large $n_i$ and $n$, is
$$4\, \breve \nu(m,\tau,c)=\frac{n^2(n-1)^2}{n_{_T}^2} 4\,\widetilde \nu(m,\tau,c)
=4\,\widetilde \nu(m,\tau,c) \, \left(\sum_{i=1}^{m+1} w_i^2\right)^{-2}
$$
since
$$\frac{n^2(n-1)^2}{n_{_T}^2}=\left(\sum_{i=1}^{m+1} n_i(n_i-1)/(n(n-1))\right)^{-2} \approx \left(\sum_{i=1}^{m+1} w_i^2\right)^{-2}$$
for large $n_i$ and $n$,
Here $\widetilde \nu(m,\tau,c)$ is as in Theorem \ref{thm:MI-asy-norm-II}.
Hence the desired result follows.
$\blacksquare$

\end{document}